\magnification1200
\font\medium=cmbx10 scaled \magstep1
\font\large=cmbx10 scaled \magstep2
\def\bn{\bigskip \noindent}
\def\sn{\smallskip \noindent}
\def\cO{{\cal O}}
\def\la{\longrightarrow}
\def\La{\longrightarrow}
\def\cF{{\cal F}}
\def\cH{{\cal H}}
\def\cL{{\cal L}}
\def\Pic{{\rm Pic}}

\def\cS{{\cal S}}
\def\cG{{\cal G}}

\def\cN{{\cal N}}
\def\cT{{\cal T}}
\def\cI{{\cal I}}
\def\F{{\cal F}}

\def\rk{\mathop{\rm rk}\nolimits}

\def\pr{\mathop{\rm pr}\nolimits}

\def\dim{\mathop{\rm dim}\nolimits}
\def\f{{\varphi}}
\def\bP{{\bf P}}
\def\bN{{\bf N}}
\def\bR{{\bf R}}
\def\bQ{{\bf Q}}

\def\bZ{{\bf Z}}
\def\Om{{\omega_{XÊ\vert C}^{-m}}}
\def\Omm{{\omega_{X \vert C}^{[m+1]}}}
\def\XY{{\omega^{-1}_{X \vert Y}}}
\def\XXY{{\omega_{X \vert Y}}}
\def\XC{{\omega_{X \vert C}^{-1}}}

\vskip .5cm 
{\centerline {{\large Threefolds with nef anticanonical bundles}}}

\vskip .5cm
{\centerline {{\medium Thomas Peternell and Fernando Serrano}}}

\vskip .5cm

\noindent {\medium Introduction}
\bn In this paper we study the global structure of projective threefolds $X$ whose anticanonical
bundle $-K_X$ is nef. In differential geometric terms this means that we can find metrics on
$-K_X = {\rm det}ÊT_X$ (where $T_X$ denotes the tangent bundle of $X$) such that the negative part
of the curvature is as small as we want. In algebraic terms nefness means that the intersection
number $-K_X \cdot C \geq 0$ for every irreducible curve $C \subset X.$ The notion of 
nefness is weaker than the requirement of a metric of semipositive curvature and is the 
appropriate notion in the context of algebraic geometry. 
\bn In [DPS93] it was proved that the Albanese map $\alpha : X \la {\rm Alb}(X)$ is a
surjective submersion if $-K_X$ carries a metric of semi-positive curvature, or, equivalently,
if $X$ carries a K\"ahler metric with semipositive Ricci curvature. It was conjectured that 
the same holds if $-K_X$ is only nef, but there are very serious difficulties with the
old proof, because the metric of semi-positive curvature has to be substituted by a sequence
of metrics whose negative parts in the curvature converge to $0.$ The conjecture splits naturally into two parts: surjectivity of $\alpha$ and smoothness.
Surjectivity was proved in dimension 3 already in [DPS93] and in general by Qi Zhang [Zh96],
using char $p.$ Our main result now proves smoothness in dimension 3:

\bn {\bf Theorem.} {\it Let $X$ be a smooth projective threefold with $-K_X$ nef. Then the 
Albanese map is a surjective submersion.}

\bn Actually much more should be true: there should be a splitting theorem
: the universal cover of $X$ should be the product of a euclidean space
${\bf C}^m$ and a simply connected compact manifold.
Again this is true if $X$ has semipositive Ricci curvature [DPS96]. 

\bn The above theorem should also be true in the K\"ahler case. Surjectivity in the threefold
K\"ahler case is proved in [DPS97], in higher dimensions it is still open. Concerning smoothness
for K\"ahler threefolds, our methods use minimal model theory, which at the moment is not
available in the non-algebraic situation.

\bn We are now describing the methods of the proof of the above theorem. First of all notice
that we may assume that $K_X$ is not nef, because otherwise $K_X$ would be numerically trivial
and then everything is clear by the decomposition theorem of Beauville-Bogomolov-Kobayashi,
see e.g. [Be83].
Since $K_X$ is not nef, we have a contraction of an extremal ray, say $\varphi: X \la Y.$ 
The Albanese map $\alpha$ factorises over $\varphi$ (of course we assume that $X$ has at least
one 1-form). If $\dim Y < 3,$ the structure of $\varphi$ is well understood and we can work
out the smoothness of $\alpha$ using the informations on $\varphi.$ So suppose that $\varphi$
is birational. It is easy to see [DPS93] that $\varphi$ has to be the blow up of a smooth 
curve $C \subset Y.$ If $-K_Y$ is nef, then we can proceed by induction on $b_2(Y).$ This is
almost always the case, but unfortunately there is one exception, namely that $C$ is rational
with normal bundle $\cO (-2) \oplus \cO(-2).$ This exception creates a lot of work; the way
to get around with this phenomenon (which a posteriori of course does not exist!), is to 
enlarge the category in which we are working. Needless to say that we have to consider threefolds
with ${\bf Q-}$factorial terminal singularities; shortly called terminal threefolds.

\bn We say that  $-K_X$ is almost nef for a terminal threefold $X,$ if $-K_X \cdot C \geq 0$
for all curves $C$ with only finitely many exceptions, and these exceptions are all rational curves.
Now in our original situation $-K_Y$ is almost nef. So we can repeat the step; if the next
contraction, say $\psi : Y \la Z$, is again birational, then $-K_Z$ will be almost nef. If
$\dim Z \geq 2,$ we can construct a contradiction: $\psi $ must be a submersion and $-K_Y$
is nef.
\sn Performing this program, i.e. repeating the process on $Z$ if necessary, we might
encounter also small contractions (contracting only finitely many curves). Then we have to
perform a flip and fortunately this situation is easy in our context, the existence of flips
being proved by Mori. Since there are no infinite sequences of flips, we will reach after a 
finite number of steps the case of a fibration $X' \la A$ and at that level the Albanese will be a
submersion. Now we still study backwards and see that we can have blown up only a finite
number of \'etale multi-sections over $A$ in case $\dim A = 1$ and that $X  = X'$ if $\dim A
= 2.$   
\bn In the last section we treat the relative situation: given a surjective map $\pi: X \la
Y$ of projective manifolds such that $-K_{X \vert Y}$ is nef, is it true, that $\pi$ is
a submersion?
\sn Our main theorem is of course the special case that $Y$ is abelian and $\pi$ the
Albanese map. We restrict ourselves again mostly to the 3-dimensional case and verify the
conjecture in several special cases. We also show that in case $\dim Y = 2,$ we may assume
that $Y$ has positive irregularity but no rational curves. However, to get around with the 
general case, we run into
the same trouble as before with the exceptional case of a blow-up of a rational curve with
normal bundle $\cO (-2) \oplus \cO (-2).$ Hopefully this difficulty can be overcome in the
near future.
  
\bn To attack the higher dimensional case however, it will certainly be necessary to
develop new methods. 

\bn We want to thank the referee for very useful comments and for pointing out some
inaccuracies.

\bn This paper being almost finished modulo linguistical efforts, the second named author
died in february 1997. Although we have never met personally, the first named author will
always remember and gratefully acknowledge the fruitful and enjoyable collaboration
by letters and electronic mail.

\vfill \eject  
{\medium 0. Preliminaries}

\bn {\bf (0.1)} Let $X$ be a normal projective threefold.
{\item {(1)} $X$ is {\it terminal} if $X$ has only terminal singularities.
\item {(2)} We will always denote numerical equivalence of divisors or curves by $\equiv.$
\item {(3)} A morphism $\f : X \la Y$ onto the normal projective threefold $Y$ is an extremal contraction
(or Mori contraction) if $-K_X$ is $\f-$ample and if the Picard numbers satisfy $\rho(X) = \rho(Y) + 1.$
\item {(4)} We let $N^1(X)$ be the vector space generated by the Cartier divisor on $X$ modulo $\equiv$ and
$N_1(X)$ the space generated by irreducible compact curves modulo $\equiv.$ 
\item {(5)} Moreover $\overline {NA}(X)
\subset N^1(X)$ is the (closure of the) ample cone, and $\overline {NE}(X) \subset N_1(X)$ is the
smallest closed cone containing all classes of irreducible curves.
\bn  In the whole paper we will freely use the results from classification theory and Mori theory 
and refer e.g. to [KMM87],[Mo87],[MP97]. The symbol $X \rightharpoonup Y$ signifies a rational morphism from $X$ to $Y.$ 

\bn {\bf (0.2)} A ruled surface is a $\bP_1-$bundle $S$ over a smooth compact curve $C.$ It is given as
$\bP(E)$ with a rank 2-bundle $E$ on $C.$ We can normalise $E$ such that $H^0(E) \ne 0$ but $H^0(E \otimes
L) = 0$ for all line bundles $L$ with negative degree. We define the invariant $e$ of $S$ by 
$e = - c_1(E).$ A section of $E$ defines a section $C_0$ of $S \la C$ with $C_0^2 = -e.$ For
details and description of $\overline {NA}(S)$ and $\overline {NE}(S)$ we refer to [Ha77,chap.V.2].
Note that $E$ is semi-stable if and only if $e \leq 0.$ 

\bn {\bf (0.3)} Let $X$ be a normal variety with singular locus $S.$ Let $X_0 = X \setminus S$ with injection
$i : X_0 \la X.$ Let $\cS$ be a reflexive sheaf of rank 1 on $X.$ Notice that $\cS$ is locally free
on $X_0.$ Let $m$ be an integer. Then we set $\cS^{[m]}Ê= i_*((\cS \vert X_0)^{\otimes m}).$  

\bn {\bf 0.4 ÊProposition}  {\it  Let $X$ be a smooth threefold, C a smooth curve and $\pi: X \la C$ a smooth
morphism and such that $-K_F$ is nef for all fibers $F$ of $\pi$ and such that $F$ is not minimal. 
Then there exists an \'etale base change $\sigma : Y = X \times_C D \la D$
induced by an \'etale map $D \la C,$ and a smooth effective divisor $S \subset Y$ such that the restriction
$ \sigma \vert S : S \la D$ yields a $\bP_1-$bundle structure on $S$, and $S \cap F$ is a 
$(-1)-$curve in $F$ for all $F.$ Hence $Y$ can be blown down along $\sigma \vert S.$}

\bn {\bf Proof.} First note that all non-minimal surfaces $F$ with $-K_F$ nef are isomorphic to the
plane $\bP_2$ blown up in at most 9 points in sufficiently general position [CP91]. Fix an ample divisor $H$
on $X.$ Pick a fiber $F$ of $\pi$ and 
take a $(-1)$-curve $E \subset F$ such that $H \cdot E$ is minimal under all $(-1)-$curves in $F.$ It follows immediately
that the normal bundle is of the form
$$ N_{E \vert X} = \cO \oplus \cO(-1).$$
By the general theory of Hilbert schemes it follows that $E$ moves algebraically in a 1-dimensional
family, i.e. there exists a projective curve $B$ and an irreducible effective divisor
$M \subset X \times B,$ flat over $B,$ such that $M \cap (X \times 0) = E,$ identifying  $X \times 0$ with $X.$ 
We let
$$ E_t = M \cap (X \times t)$$
and shall identify $X$ with $X \times t.$ 

\bn {\it Claim 1} \ ÊEvery $E_t$ is a $(-1)-$curve in some fiber $F'$ of $\pi.$

\noindent It is clear that $E_t$ has to be a Cartier divisor in some fiber $F'$ (consider the
deformations of the line bundle $\cO_F(E)$). In particular every $E_t$ is
Gorenstein and Cohen-Macaulay and does not have embedded points. Observe next that
$$ -K_X \cdot E_t = -K_X \cdot E = 1.$$
If $-K_F$ is ample for every $F$, then we deduce that $E_t$ is irreducible and reduced and by flatness that
$E_t \simeq \bP_1.$ Hence $E_t$ is a $(-1)$curve in $F'.$ If $-K_F$ is merely nef, we need to be more
careful. Assume that some $E_{t_0}$ is reducible. Write
$$ E_{t_0} = \sum a_i C_i $$
with irreducible curves $C_i.$ Since $-K_{F'}$ is nef, we conclude (after renumbering possibly) that
$$a_0 = 1, -K_{F'}Ê\cdot C_0 = 1$$ 
and that 
$$-K_{F'}Ê\cdot C_i = 0, i \geq 1.$$ 
We claim that $H^1(\cO_{E_{t_0}}) = 0$ and therefore that all $C_i$ are smooth rational curves.
One is tempted to argue by flatness, however it is not completely clear that $h^0(\cO_{E_{t_0}}) = 1,$
since $E_{t_0}$ might not be reduced. So we argue as follows. Consider the exact sequence
$$ H^1(\cO_{F'}) \la H^1(\cO_{E_{t_0}}) \la H^2(\cI_{E_{t_0}}).$$
Since $F'$ is rational, it suffices to see $$H^2(\cI_{E_{t_0}}) = 0. \eqno (*)$$ Note that
$$H^2(\cI_{E_{t_0}}) \simeq H^0(F',\cO(E_{t_0}) \otimes \omega_{F'}).$$
Now $F'$ is realised as blow-up of $\bP_2$ in 9 points. Therefore it makes sense to speak of a
general line in $F'.$ Take such a general line $l$ in $F'.$ It can be deformed to a general line
$l_s$ in a neighboring $F_s.$ Now for general s we have
$$ (K_{F_s} + E_s)Ê\cdot l_s < 0$$
where $E_s$ is one of the $(-1)-$curves in our family sitting in $F_s.$ Therefore
$$ (K_{F'} + E_{t_0}) \cdot l < 0$$ proving (*). We conclude in particular that
$C_0 $ is a $(-1)-$curve in $F'$ and the $C_i, i \geq 1,$ are $(-2)-$curves.
We claim that $B$ is smooth at $t_0.$ For this we need to know
$$ h^1(N_{E_{t_0}Ê\vert X}) = 0.$$
This comes down to $$ h^1(N_{E_{t_0}Ê\vert F'}) = 0$$ 
since $h^1(\cO_{E_{t_0}}) = 0.$ By $H^q(\cO_{F'}) = 0, q \geq 1,$ we must prove
$$ h^1(\cO_{F'}(E_{t_0})) = 0.$$
If this would not be true, then by $\chi(\cO_{F'}(E_{t_0})) = 1$, $E_{t_0}$ would move inside $F'.$ Any deformation of $E_{t_0}$ 
must have however the same type of decomposition, so that necessarily some of the $C_i$ 
would have to move in $F'$ which is absurd.

\noindent  Now we look at the deformations of $C_0$ and obtain a family 
$(C_s)_{s \in A.}$ For a small neighborhood $\Delta \subset B$ of $t_0$ the curve $E_t $ is in $\pi^{-1}(t)$
(strictly speaking there is a canonical map $f: B \la C,$ and $f\vert \Delta $ is an isomorphism, so that
we can identify $t$ and $f(t)$ for small $t).$ In the same way, $C_t \subset \pi^{-1}(t).$ 
Therefore we can consider the (non-effective) family of cycles $(E_t - C_t)_{t \in \Delta}$ so that
$$ E_{t_0} - C_0 = \sum_{i \geq 1} a_i E_i.$$ 
By the choice of $E_t,$ $H \cdot E_t$ is minimal for general $t,$ therefore $H \cdot E_t \leq H \cdot C_t$ and
we conclude 
$$ H \cdot \sum_{i \geq 1} a_iÊE_i = 0$$
and therefore $a_i = 0$ for $i \geq 1$ so that $E_0$ is irreducible and reduced. 

\noindent This proves Claim 1.

\bn {\it Claim 2} Let $Z = \pr_1(M) \subset X.$ Then $Z \cap F'$ is a reduced union of $(-1)-$ curves
and the number is independent of $F'$. 

\noindent In fact, the first part (reducedness) is immediate from Claim 1 (if a $(-1)-$curve $E$ in a fiber appears with multiplicity
$m \geq 2$ in $Z \cap F,$ then $E$ could be deformed itself to the neighboring fibers. This 
contradicts clearly the smoothness of
$B)$. The independence of the number follows also 
from the smoothness of $B.$

\noindent In other words, Claim 2 says that $f: B \la C$ is \'etale. So set $D = B, \  Y = X \times_C D$ and
define $S$ to be the irreducible component of $  Z \times_C D$ mapping onto $Z.$ 

\bn {\bf (0.5) Remark.} In (0.4) we used the nefness assumption for $-K_X$ only to make sure
that (-1)-curves in fibers can only be deformed into (-1)-curves in fibers. If we know this
for some other reason, then the conclusion of (0.4) remains true.
 
\bn The next proposition should be well-known and hold in more generality; however we could not 
find a reference, so
we include the short proof. 
\bn {\bf (0.6) Proposition} {\it Let $X$ be a terminal ${\bf Q}$-factorial threefold and et $\f: X \la Y$ be the contraction of an extremal ray.
Assume that $\f$ contracts a divisor $E$ to a curve $C.$ Assume that $Y$ is smooth and $C$
is locally a complete intersection. Then $\f$ is the blow-up of $C$ in $Y.$}

\bn {\bf Proof.} Let $N$ be the singular locus of $C$ and $\tilde N = \varphi^{-1}(N).$
Then $\tilde N$ is purely 1-dimensional or empty since $E$ is irreducible. 
Let $\pi : X' \la Y$ be the blow-up of $Y$ along $C$ with exceptional set $E'.$ Since
$C$ is locally a complete intersection, we have $E' = \bP ({\cal N}^*_C).$ Thus 
$N' = \pi^*(N)$ is purely 1-dimensional or empty, too. Since $\varphi$ is generically the blow-up
of $C,$ we have an isomorphism $\tilde X \setminus \tilde N \La X' \setminus N'.$
We next observe that $X'$ is normal. Locally (in $Y$) we have $X' \subset Y \times \bP_1,$ since $Y$ is smooth.
Hence $X'$ is Cohen-Macaulay. On the other hand, $$\dim {\rm Sing}(X') \leq 1.$$
In fact, up to a finite set, ${\rm Sing}(X') \subset \pi^{-1}({\rm Sing}(C)).$
Now all non-trivial fibers of $\pi$ are (smooth rational) curves, hence
$\dim {\rm Sing}(X') \leq 1.$ 
\sn Putting things together, $X'$ is normal.
By [Ko89,2.1.13] we have $X \simeq X'$ unless $N'$ has an irreducible contractible component which
is of course absurd.     

\vfill \eject

{\medium 1. Fiber spaces}
\bn 
{\bf 1.1 Definition }ÊÊLet $X$ be a normal projective variety and $L$ a line bundle on $X.$
Then $L$ is called {\it almost nef}, if there are at most finitely many rational curves $C_i,
1 \leq i \leq r,$ such that $L \cdot C \geq 0$ for all curves $C \ne C_i.$ 

\bn {\bf 1.2 Proposition}  {\it Let $X$ be a terminal n-fold with $-K_X$ almost nef. Then $\kappa (X) \leq 0.$ 
Moreover the following three statements are equivalent. 
{\item {(a)} $\kappa (X) = 0$
\item {(b)} $K_X \equiv 0$
\item {(c)}Ê$K_X$ is nef}}

\bn  {\bf Proof.} The first assertion is clear. 
\sn If $\kappa (X) = 0$ and if $K_X \not \equiv 0,$ then there exists a non-zero $D \in
\vert mK_X \vert$ for some positive $m.$  Hence $-K_X$ cannot be almost nef.
\sn If $K_X $ is nef, then $K_X \cdot C = 0$ for all but finitely many curves. In particular
we have $K_X \cdot H_1 \cdot ... \cdot H_{n-1} = 0$ for all ample $H_i$ on $X.$ Therefore
$K_X \equiv 0,$ see e.g. [Pe93,6.5].

\bn {\bf 1.3 Proposition} {\it Let $X$ be a terminal ${\bf Q}$-factorial 3-fold with $-K_X$ almost nef. Assume
that there is an extremal contraction $\varphi: X \la C$ to the elliptic curve $C.$
Then $X$ is smooth, $\varphi $ is a submersion and $-K_X$ is nef. }

\bn {\bf Proof.} All rational curves in $X$ are contracted by $\f,$ moreover all rational
curves are homologous up to multiples. Hence $-K_X$ must be nef.
\sn  
\sn (A) First note that for all positive $m$ the sheaf 
 $$V_m = \f_*(-mK_{X\vert C}) = \f_*(\omega_{X\vert C}Ê^{[-m]}) $$  is a 
vector bundle since it is torsion free and $\dim C = 1.$ 
Now let us consider only those $m \in \bN$ such that $mK_X$ is Cartier. 
Then we have 
$$ c_1(V_m) \leq 0. \leqno (*)$$
For the proof of (*) we first compute (using the relative version of Kawamata-Viehweg, recall that $\omega_{X \vert C}^{-1}$
is $\varphi-$ample) 
$$ \chi(\Om) = \chi(\f_*(\Om)) = \chi(V_m) = c_1(V_m) $$
by Riemann-Roch on $C.$ Next we compute
$\chi(\Om) $
on $X.$ The first step is to apply Riemann-Roch to obtain 
$$ \chi(\Omm) = \chi(\omega_X^{[m+1]}) = (1-2(m+1)) \chi(X,\cO_X) + A ,$$
where $A \geq 0$ (and $A = 0$ if and only if $X$ is Gorenstein); see [Re87],[Fl87] for the singular
Riemann-Roch version needed here. Note that we have used $K_X^3 = 0;$ in fact, if $K_X^3 < 0,$ 
then $-K_X$ would be big and nef, hence $q(X) = 0$ (see [KoMiMo92,3.11]. Since
$$ \chi(\cO_X) = \chi(\cO_C) = 0,$$
we get $$ \chi(\Omm) \geq 0. \eqno (1)$$
Since $mK_X$ is Cartier, we have
$$ \Omm = \omega^m_{X \vert C} \otimes \omega_{X \vert C} = \omega^m_{X \vert C}Ê\otimes \omega_X,$$
hence
$$ \chi(X,\Omm) = - \chi(X,\omega^{-m}_{X \vert C}) $$
by Serre duality. Thus $\chi(\Om) \leq 0$ and we conclude $c_1(V_m) \leq 0.$

\sn (B) We claim that $V = V_m$ is nef. In case $X$ is smooth and $\f$ a submersion this is just [DPS94,
3.21], applying (3.21) to $L  = \omega^{-(m+1)}_{X \vert C}.$ The proof of (3.21)  remains 
valid in our situation if $\f$ is only flat (which is true since $\dim C = 1,$  but $X$ is still
assumed to be smooth). If $X$ is singular, we argue as follows. Let $\pi : \hat X \la X$ be
a desingularisation and let $\hat \f : \hat X \la C$ denote the induced map. Let 
$$L_m = \pi^*(\omega_{X \vert C}^{[-(m+1)]}) / {\rm torsion};$$
then $L_m$ is locally free, at least if $\pi$ is chosen suitably (see e.g. [GR70]). 
At the same time we can achieve that
$$ \pi^*(\omega_{X \vert C}^{[-(m+1)]})^{\otimes m}Ê/ {\rm torsion} = 
\pi^*(\omega_{X \vert C}^{[-m(m+1)]}) / {\rm torsion}$$
is locally free. Then it is immediately checked that $mL_m = \pi^*(\omega_{X \vert C}^{-m})^{m+1},$
therefore $mL_m$ is nef, and so does $L_m.$   
 By the flat version of [DPS94,3.21] the bundle
$$ \hat \f _* (\omega_{\hat X \vert C}Ê\otimes L_m) $$
is nef. Now
$$  \pi_*(\omega_{\hat X \vert C}Ê\otimes L_m) \subset (\omega_{X \vert C} \otimes 
\omega_{X \vert C}^{[-(m+1)]})^{**} = \Om,$$
since the first sheaf is torsion free, the second is reflexive and both coincide outside
a finite set. Therefore
$$ \hat \f _* (\omega_{\hat X \vert C}Ê\otimes L_m) \subset \f _*(\Om) = V_m,$$
and the inclusion is an isomorphism generically. Thus $V_m$ is nef. 
Since $c_1(V_m) \leq 0,$ we conclude that $V_m$ is numerically flat, i.e. both $V_m$ and
$V_m^*$ are nef (see [DPS94]), in particular $c_1(V_m) = 0.$    

\sn By (B) we conclude
$$ \chi(X,\omega_X^{[m+1]}) = - \chi(X,\omega_X^{-m}) = \chi(V_m) = 0.$$
Therefore our reasoning in (A) proves that $X$ is Gorenstein. 
 
\sn (C) If $m \gg 0,$ we have an embedding
$$  i: X \hookrightarrow \bP (V),$$
since $-mK_{X \vert C}$ is $\f-$very ample. Let $r = \rk V$ and $\cO_X(1) = i^*(\cO_{\bP (V)}(1)).$
Then by construction
$$ -mK_X = \cO_X(1) \otimes \f^*(L)$$ 
with some line bundle $L$ on $C.$ 
We claim that 
$$c_1(L) = 0.$$
To verify this, first notice from $-mK_X = \cO_X(1) \otimes \varphi^*(L)$ that 
$$  V = V_m = \varphi_*(\cO_X(1)) \otimes L. \eqno (+)$$
Now consider the exact sequence
$$  0 \la \cI_X \otimes \cO_{\bP(V)}(1) \la \cO_{\bP(V)}(1) \la \cO_X(1) \la 0 $$
and apply $\pi_*$ to obtain 
$$ 0 \la \pi_*(\cI_X \otimes \cO_{\bP(V)}(1)) \la V \la \varphi_*(\cO_X(1)) 
\la R^1\pi_*(\cI_X \otimes \cO_{\bP(V)}(1)) \la 0 \eqno (++).$$
We check that 
$$ R^1\pi_*(\cI_X \otimes \cO_{\bP(V)}(1)) = 0.$$
In fact, this sheaf is $0$ generically, since for general $c \in C$, the embedding $X_c = \varphi^{-1}(c) \subset \pi^{-1}(c)$
is defined by $H^0(X_c,-mK_{X_c})$ which implies 
$$H^1(X_c, \cI_{X_c} (1)) = 0.$$
Since however $\varphi_*(-mK_X) $ is locally free and $R^q\varphi_*(-mK_X) = 0$ for $q > 0,$ standard semi-continuity
theorems (notice that $\varphi$ is flat!) imply that $h^0(X_c,-mK_{X_c})$  is constant. Since 
$$ H^1(X_c,-mK_{X_c}) = 0,$$
as we check easily, we obtain 
$$ H^1(X_c,\cI_{X_c}(1)) = 0$$
for all $c \in C,$ hence $R^1\pi_*(\cI_X \otimes \cO_{\bP(V)}(1)) = 0.$
Since ${\rm rk}V = {\rm rk}(\varphi_*(\cO_X(1))$ by (+), we conclude 
$$ V \simeq \varphi_*(\cO_X(1))$$ by (++) and the $R^1-$vanishing. Again from (+) we finally obtain
$$ c_1(L) = 0.$$ 
 
Using
$$ K_{\bP (V)}Ê= \cO_{\bP (V)}(-r) \otimes \pi^*(\det V),$$
we obtain by the adjunction formula
$$ \cO_X(rm-1) = \f^*((\det V)^m \otimes L) \otimes (\det N_X)^m.  \leqno (**)$$
Now it is well-known that $\bP (V) $ is almost homogeneous (and the tangent bundle
$T_{\bP (V)}$ is nef) (cp. [CP91]), i.e. the holomorphic vector fields generate $T_{\bP (V)}$
outside some proper analytic set $S \subset \bP (V).$ 

\sn (C1) We first treat the case that $ X \not \subset S.$ Assume that 
$\f$ is not a submersion. This means that the sheaf of relative K\"ahler differentials 
$\Omega^1_{X \vert C}$ is not locally free of rank 2. Note that once we know that $\Omega^1_{X \vert C}$
is locally free, then automatically $X$ must be smooth. We are first going to show that under our assumption 
$$ h^0(\cN_{X \vert \bP(V)}) > \rk \cN. \leqno (2)$$ 
Here $\cN$ denotes the normal sheaf of $X \subset \bP (V),$ the dual of $\cI / \cI ^2.$
Let $\omega$ be the pull-back of a non-zero 1-form from $C.$ From the exact sequence 
$$ 0 \la \f^*(\Omega^1_C) \la \Omega^1_X \la \Omega^1_{X \vert C} \la 0$$ 
we see that $\omega$ has zeroes exactly at some of the singularities of $X$ and at the smooth points
of $X$ where
$\f$ is not a submersion. Consider the exact sequence of tangent sheaves
$$ 0 \la \cT_X \la \cT_{\bP (V)}Ê\vert X \la \cN_X. \leqno (S) $$
This sequence shows that $\cN_X$ is generated by global sections outside the set
$\tilde S = S \cup {\rm Sing}ÊX.$ If $$h^0( \cN) = {\rk} {\cN} ,$$
then by (S) also $\cT_X$ would be generically generated. Hence we can find $v \in H^0(\cT_X)$
such that $\omega (v) \ne 0,$ so that $\omega (v)$ is a non-zero constant holomorphic 
function and $\omega$ has no zeroes. Therefore $\f$ can fail to be a submersion at most at the singularities of
$X$, in particular inequality (2) holds already for $X$ smooth. In the remaining case we argue as follows.
Since $\cT_X$ is generically generated, $X$ is almost homogeneous with respect to ${\rm Aut}^{o}(X),$  i.e.
the automorphisms act with an open orbit. Every $x \in {\rm Sing}(X)$ must be a fixed point. Hence the fiber
of $\f$ containing $x$ is invariant under the action and consequently the induced action on $C$ has a 
fixed point. $C$ being elliptic, the action on $C$ is trivial, but then $X$ cannot be almost homogeneous.
Of course this argument can also be used in the case $X$ smooth. 
\sn Now (2) is proved. In particular, $\cN_X$
being generically spanned, we have
$$ h^0(\det \cN_X )  \geq 2.$$
By (**) we conclude the existence of some $n_0 \in \bN$ and a line bundle $\cG_0 \in 
\Pic^0(C)$ such that 
$$ h^0(\omega_X^{-n_0}Ê\otimes \f^*(\cG_0)) \geq 2.$$
Note that necessarily $n_0 K_X$ is Cartier.   
We claim:
\sn (***) there is some $n_1\in \bN$ and a $\cG_1 \in \Pic^0(C)$ such that the base
locus $B_1$ of the linear system $\vert  \omega_X^{-n_1}Ê\otimes \f^*(\cG_1) \vert $
has dimension $\leq 1.$ 
\sn {\it Proof of (***):} If already the base locus $B_0$ of our linear system
$\vert  \omega^{-n_0}_X \otimes \f^*(\cG_0) \vert $ has dimension $\leq 1,$ then we are done; so
assume that $\dim B_0 = 2$. Let $\tilde B_0$ be the
2-dimensional part (with appropriate multiplicities). Let
$$ M = \omega_X^{-n_0}Ê\otimes \f^*(\cG_0) \otimes \cO_X (- \tilde B_0).$$
Then the base locus of $\vert M \vert $ has dimension at most $1.$ We can write
$$ \cO_X(\tilde B_0) = \omega^{-\mu}_X \otimes \f^*(H),$$
note that $\tilde B_0$ is Cartier and that $\mu$ is a non-negative rational number and $H$ is
$\bQ-$Cartier on $C.$ Now choose $k \in \bN$ such that $k(n_0 + \mu) = \rho m$ for some
positive integer $\rho $ where $mK_X$ is Cartier and let $n_1 = \rho m.$ Then
we consider $kM$ instead of $M,$ of course the base locus of $\vert kM \vert $ still has  dimension
at most $1.$ We have
$$ kM = \omega_X^{-n_1} \otimes \f^*( \cG_0^k \otimes H^{-k}).$$ 
If $H \equiv 0,$ we are done, so assume $H \not \equiv 0.$ Since
$$ 0 \ne H^0(X,M^k) = H^0(C,V_{n_1} \otimes H^{-k} \otimes \cG_0^k),$$ 
the numerical flatness of $V_{n_1}$ forces $\deg H < 0.$ But then, going back to the
decomposition of $\tilde B_0,$ we would have a
section of $-\mu K_X$ vanishing on some fibers of $\f$ which gives a section of $V_{\mu}$ 
with zeroes, contradicting the flatness of $V_{\mu}.$ So $H \equiv 0.$ 
This proves (***). 
\sn  Let
$$f : X \rightharpoonup Y $$ be the map associated to the linear system
$ \vert \omega_X^{-n_1}Ê\otimes \f^*(\cG_1) \vert.$ Since $-K_X$ is not big, we have
$\dim Y = 1$ or $\dim Y = 2.$  Let $F$ be a general fiber of $f.$ Note first that in case $\dim Y = 1,$
the map $f$ cannot be holomorphic, i.e. $B_1 \ne \emptyset,$ because otherwise $K_X^2 = 0,$ which is
impossible, $\f$ being a del Pezzo fibration. We next treat the case that $f$ is holomorphic in case
$\dim Y = 2,$ or, more generally, that $F \cap B_1 = \emptyset.$ Then
$$ K_F = K_X \vert F \equiv 0.$$ 
Therefore $F$ is an elliptic curve. 
Moreover $$c_1(\cO_{\bP (V)}(1)) \cdot F = 0.$$ 
Using the tangent bundle sequence and the generic spannedness of $\cT_{\bP (V)},$ we see
immediately that $\cN_{F\vert \bP (V)} = \cO_F^{\oplus N}.$ 
Now the relative tangent bundle sequence for $ \pi :\bP (V) \la C$ together with the relative Euler
sequence imply that
$$ h^*(V) = \cO_F^N,$$
where $h = \pi \vert F \la C$ is the \'etale covering of $F$ over $C.$ Hence after the base change
$F \la C$ the space $\bP (V)$ becomes a product. It follows in particular that $f$ must be
an elliptic bundle and that $\f$ is smooth.
\sn So we are reduced to the case that $B_1 \ne \emptyset$. Then we even have $\dim B_1 = 1,$
otherwise we could pass to $m(-n_1 K_X + \f^*\cG_1))$ to obtain base point freeness. Let $B \subset
B_1$ be the 1-dimensional part of $B_1.$
\sn (a) We start with the case $\dim Y = 1.$ First note that
$$ F \equiv - \rho K_X$$
with some positive rational number $\rho.$ Take another general fiber $F'$ and consider the nef line
bundle $F' \vert F$ (strictly speaking we should take $\lambda$ such that $\lambda F'$ is Cartier and
consider $\lambda F' \vert F).$ We write (on $F$)
$$ F' \vert F = B + M,$$
$M$ the movable part. Decomposing $B = \sum b_i B^{i}, $ we deduce
from $K_X^3 = 0$ that
$$ -K_X \cdot \sum b_i B^{i}Ê+ M = 0.$$
Since $-K_X$ is nef, we conclude $-K_X \cdot B^{i} = -K_X \cdot M = $ for all $i.$ Therefore all $B^{i}$
and $M$ are homologous, i.e. contained in the  half ray 
$$R = \{ÊZ \in \overline {NE}(X) \vert Z \cdot K_X  = 0 \}$$ 
inside the 2-dimensional cone $\overline {NE}(X).$ 
There is a slight difficulty that $M$ and $B$ a priori might not be $\bQ-$Cartier in $F.$ To
circumvent this, choose a desingularisation $\sigma: \hat F \la F.$ Let $\hat M$ be the strict
transform of $M$ in $\hat F.$ Choose $\hat B_j \subset \hat F$ such that $\sigma (\hat B_j) \subset
B^{i(j)}$ and such that there is an equation
$$ \sigma^*(F' \vert F) = \hat M + \sum \hat b_j \hat B_j + E, \leqno (3)$$
where $E$ is effective and contained in the exceptional locus for $\sigma$ (including the non-normal
part). $\hat M$ being irreducible and movable (for general choice of $M$), we have $\hat M^2 \geq 0.$
If $\hat M^2 > 0,$ then $\hat M$ would be big, so $\sigma(F' \vert F)$ would be big contradicting the
nefness of $F'$ together with $F'^2 \cdot F = 0.$ Hence $\hat M^2 = 0.$ Thus $\hat M$ is base point free
and defines a map
$$ \hat \lambda : \hat F \la B_F$$
to a curve $B_F.$ Now notice 
$$ \sigma^*(F') \cdot \hat M = F' \cdot M = 0$$
(use $F'^2 \cdot F = 0$ and the nefness of $F' \vert F$).
Therefore $\sigma^*(F') \cdot l = 0,$ with $l$ a fiber of $\hat \lambda.$ Consequently all $\hat B_j$
and all components of $E$ must be contained in fibers of $\hat \lambda$ (just dot (3) with $l)$.
It follows that $M$ is Cartier on $F$ (and so does $B)$ and its sections define a morphism
$$\lambda_F : F \la B'_F$$
to a curve $B'_F$ (with a natural map $B_F \la B'_F)$. 
Since $M \cdot B^{i} = 0$ (in $F$), all $B^{i}$ are contracted by $\lambda_F$ and hence the general
fiber $G$ of $\lambda_F$ does not meet $B_1.$ We may assume $G$ connected.  Since
$$K_G = K_F \vert G \equiv (1 - \rho)K_X \vert G,$$
we have $K_G \equiv 0$ and either $G$ is  smooth elliptic or a singular rational curve. This second
alternative cannot occur: since $\dim \f(G) = 1$ by virtue of $K_X \cdot G = 0,$ the curve $G$ surjects
to the elliptic curve $C.$ Hence $G$ is a smooth elliptic curve. Now we argue as in the case $\dim Y = 2$
and $f$ holomorphic and obtain a contradiction.
\sn (b) The case $\dim Y = 2$ with $\dim B_1 = 1$ is essentially the same. We choose 
$$ D,D' \in \vert  -n_1 K_X + \f^*(\cG) \vert $$
general, subsitute $F$ by $D$ and $F'$ by $D'$ and repeat the arguments of (a).          
This finishes the case (C1).
\sn (C2) We still must deal with the case $X \subset S.$ 
The structure of $S$ is however very easy. Choose $\cH \in \Pic^0(C),$ such that, putting
$\tilde V = V \otimes \cH,$ the dimension $h^0(\tilde V)$ gets maximal. Write $\tilde V$
as the following extension
$$ 0 \la \cO_C^p \la \tilde V \la V' \la 0,$$
such that $h^0(V') = 0.$ Then the exceptional orbit $S$ is of the form $S = \bP (V') \subset
\bP (\tilde V) = \bP (V).$ Now we substitute $V$ by $V'$ and run the old argument. 

\bn We  proceed with investigating conic bundles over possibly singular surfaces.

\bn {\bf 1.4 Lemma} {\it Let $Y$ be a normal projective surface with only rational singularities.
Assume that $-K_Y$ is almost nef and that $q(Y) \geq 1.$ Then $Y$ is either a $\bP_1-$bundle over an 
elliptic curve, an abelian surface or a hyperelliptic surface; in particular $Y$ is smooth, 
and $-K_Y$ is nef. }

\bn {\bf Proof.} Rational singularities are automatically ${\bQ-}$Gorenstein, hence the assumption
"$-K_Y$ nef" makes sense. Let $\pi: \hat Y \la Y$ be a desingularisation. Since
$$ -K_{\hat Y} = \pi^*(-K_Y) + A $$
with $A$ effective (possibly 0), $-K_{\hat Y}$ is almost nef. Let $\sigma: \hat Y \la Y_m$ be a map
to the minimal model. Then
$$ -K_{\hat Y} = \sigma^*(-K_{Y_m}) - E $$
with $E$ effective, hence $-K_{Y_m}$ is almost nef. 

If $\kappa (Y_m) \geq 0,$ we conclude that $K_{Y_m} \equiv 0,$ so that $Y_m$ is abelian or
hyperelliptic (by the existence of a 1-form); moreover that $\hat Y = Y = Y_m$ by almost nefness
of the corresponding canonical bundles. 

Hence we shall assume $\kappa (Y_m) = -\infty$ from
now on.  
$Y_m$ being a $\bP_1-$bundle over a curve $C$  of
genus $g(C) \geq 1,$ it is clear that $-K_{Y_m}$ is nef, hence $C$ is an elliptic curve. 
It remains to prove the following
\sn (*)  if $\lambda: Y' \la Y_m$ is the blow-up of the point $p \in Y_m,$ then $-K_{Y'}$ is  not
almost nef. 
\sn Given (*), we conclude that $\hat Y = Y_m,$ and since $Y$ has only rational singularities, it
follows $Y = \hat Y = Y_m.$ 
\sn For the proof of (*), we first note that $-K_{Y'}$ must be nef if it is almost nef. In fact, otherwise 
there is a rational curve $C$ with $K_{Y'}Ê\cdot C > 0.$ Since $C$ does not move, it can only be the
exceptional curve for $\lambda$ or the strict transform of the ruling line containing $p.$ But in both
cases $K_{Y'}Ê\cdot C = -1.$ Hence $-K_{Y'}$ is nef. On the other hand $ K_{Y'}^2 = -1,$ contradiction.
This proves (*) and finishes the proof of (1.4).

\bn {\bf (1.5)} Let $\f : X \la W$ be an extremal contraction of the terminal ${\bf Q}$-factorial threefold $X$ to the surface
$W.$ It is well-known and easy to prove that $\f$ is equidimensional (since $\rho (X) =  \rho (Y) + 1.$)
The surface $W$ has only quotient singularities, i.e. $(W,0)$ is log terminal, in particular $W$ has only
rational singularities (see [KoMiMo92]). Let $$S = {\rm Sing}(X); S' = \f (S)$$ and $$W_0 = W \setminus S',
X_0 = \f ^{-1}(W_0).$$ 
Then $\f _0 : X_0 \la W_0$ is a usual conic bundle and $W_0$ is smooth. Let $\Delta _0$ denote its
discriminant locus and put $\Delta = \overline {\Delta _0} \subset W.$  

\bn {\bf 1.6 Lemma}  {\it Assume the situation of (1.5). If $-K_X$ is almost nef, then $-(4K_W + \Delta)$
is almost nef.}

\bn {\bf Proof.} Note that $W$ is $\bQ-$factorial since it has only rational singularities. The arguments in
[Mi83,4.11] show that 
$$ \f_*(K_X^2) = -(4K_W + \Delta) $$
in $N^1(W),$ since this has only to be checked on curves which are very ample divisors on $W$ (and therefore
may be assumed not to pass through $S').$ Hence our claim is clear: if $$-(4K_W + \Delta) \cdot C < 0,$$
then $ K_X^2 \cdot \f^{-1}(C) \cdot C < 0,$ hence $-K_X \vert \f^{-1}(C)$ cannot be nef, hence 
$\f^{-1}(C)$ contains one of the finitely many rational curves $C'$ with $K_X \cdot C' > 0$ so that
$C = \f(C').$ 

\bn {\bf 1.7 Proposition} {\it Let $X$ be a terminal ${\bf Q}$-factorial threefold with $-K_X$ almost nef. Assume $q(X) = 1$
and let $\alpha : X \la C$ be the Albanese map to the elliptic curve $C.$ Let $\f : X \la W$ be an extremal
contraction to the surface $W.$ Then $X$ is smooth and $\alpha$ is a submersion. Moreover $W$ is a hyperelliptic
surface or a $\bP_1-$bundle over $C$ with $-K_W$ nef. }

\bn Note that we do not claim here that $-K_X$ is nef; we will address to this point in (1.8).

\bn {\bf Proof.} We shall use the notations of (1.5). If $\kappa(\hat W) \geq 0, \hat W$ 
a desingularisation, then, $-(4K_W + \Delta)$ being almost nef, $-K_W$ is the sum of an almost
nef and an effective divisor which includes $\Delta.$ Passing to $\hat W$ and using the effectiveness
of $K_{\hat W},$ it follows immediately $\Delta = 0.$ Hence $W$ is hyperelliptic by (1.4). 
But then by a base change we pass to the case ${\rm alb}(X) = 2, \dim W = 2$ treated in (1.9) and
(1.10). However it is also possible to make the following arguments work also in the hyperelliptic
case. From now we will assume $\kappa (\hat W) = - \infty.$  
\sn (A)  First we consider the case $\Delta = 0.$ By (1.6) $-K_W$ is almost nef, hence $W$ is smooth by
(1.4). We claim that
$$ X_0 = \bP (E_0) $$
with an algebraic vector bundle $E_0$ on $W_0.$ First we show that $E_0$ exists as a holomorphic bundle.
The obstruction for the $\bP_1-$bundle $X_0 \la W_0$ to be of the form $\bP (E_0)$ is a {\it torsion} element
$$ P \in H^2(W_0, \cO^*) $$
(see e.g. [El82]). From the exponential sequence we see
$$ H^2(W_0,\cO^*) \simeq H^3(W_0,\bZ),$$
if $S' \ne \emptyset.$ Assuming $S' \ne \emptyset$ for the moment, we check easily via Mayer-Vietoris that $H^3(W_0,\bZ)$ is
torsion free.  Hence $P = 0.$ If $S' = \emptyset$  then $X$ is smooth and $\f$ is a $\bP_1-$bundle so that
$\alpha$ is a submersion. Hence we will assume that $S' \ne \emptyset,$ i.e. that $X$ is singular. 
\sn Now we have $X_0 = \bP(E_0)$ analytically. Therefore $-K_{X_0 \vert W_0} = \cO_{\bP (V)}(2)$ analytically with some 
rank 2-vector bundle $V.$ We may assume $V = E_0.$ Of course $-K_{X_0 \vert W_0}$ is algebraic; we want to show that
$E_0$ is algebraic, i.e. $\cO_{\bPÊ(E_0)}(1)$ is algebraic. In fact, taking roots, there is a 2:1 Galois cover 
$g: \tilde X_0 \la X_0$ and an algebraic line bundle $\cL $ on $\tilde X_0$ such that
$$ g^*(-K_{X_0 \vert W_0})Ê= \cL^2.$$
So $g^*(-K_{X_0 \vert W_0})$ is algebraic and so does $g_*g^*(-K_{X_0 \vert W_0}) \simeq \cO_{X_0}Ê\oplus 
-K_{X_0 \vert W_0}.$                 
So $E_0$ can be taken to be an algebraic vector bundle. Thus $E_0$ has a coherent extension
to $W.$ The bidual of this extension is reflexive, hence locally free, $W$ being a smooth surface. Thus
$E_0$ has a vector bundle extension $E.$ Let $\tilde X = \bP (E).$ Then $\tilde X$ and $X$ coincide
outside finitely many curves. Thus $\tilde X \simeq X$ by [Ko89,2.1.13]. Hence $\f$ and therefore $\alpha$
is a submersion. 
\sn (B) Now let $\Delta \ne 0.$  
\sn By (1.6), $-(4K_W + \Delta)$ is almost nef. It follows already that $-K_W$ is almost nef except possibly
for the case that there might be an irrational curve $B \subset \Delta$ with $K_W \cdot B > 0.$ 
\sn If $-K_W$ is almost nef, then by (1.4) $W$ is smooth, in fact a $\bP_1-$bundle over $C$ with $-K_W$
nef [CP91].
We therefore shall prove now that $-K_W$ is almost nef. Assume to the contrary that there
is an irrational curve $B$ such that $$K_W \cdot B > 0.$$
We have already seen that necessarily $B \subset  \Delta.$ Note that $W$ is ${\bQ-}$factorial
since $W$ has only rational singularities. In particular $K_W$ and $B$ are $\bQ-$Cartier.
We claim that
$$ K_W + B \cdot B < 0.$$
In fact, since $-4(K_W + \Delta)$ is almost nef and $B$ irrational, we have 
$$ -4 (K_W + \Delta) \cdot B \geq 0,$$
so that $\Delta \cdot B \leq -4 (K_W \cdot B).$ Consequently
$$ B^2 \leq \Delta \cdot B \leq -4(K_W \cdot B) < - (K_W \cdot B).$$ 
This proves the claim. Now let $\mu : \tilde B \la B$ be the normalisation. Choose
$m$ positive such that $m(K_X + B)$ is Cartier. Then by the subadjunction lemma (see
[KMM87,5-1-9]), there is a canonical injection
$$ \omega^m_{\tilde B}Ê\la \mu^*(\cO_B(m(K_X + B))).$$
Hence ${\rm deg}K_{\tilde B}Ê< 0$ and $B$ is rational, contradiction. So $-K_W$ is almost
nef.
\sn (C) Now we know that $W$ is a $\bP_1-$bundle over $C$ with $-K_W$ nef. Hence $e(W) = 0, -1.$ 
Moreover $-(4K_W + \Delta)$ is nef. However $X$ maybe still be singular and $-K_X$ only
almost nef. First let us see that $X$ is Gorenstein and $\f$ really a conic bundle. We shall
use the notations from (A). The sheaf $$ \F = \f _{0*}(\omega_X^*)$$ is torsion free and locally free on
$W_0.$ We claim that $\cF$ is actually reflexive. In fact, take $x \in W \setminus W_0$, let
$U \subset W$ be an open neighborhood of $x$ and take $s \in H^0(U \setminus \{ x\},\cF).$ We need to prove
that $s$ extends to $U.$ Consider $s$ as an element of $H^{0}(\f^{-1}(U \setminus \{ x \},\omega_X^*).$ 
Since $\dim \f^{-1}(x) = 1$ and since $\omega_X^* = \cO_X(-K_X)$ is reflexive, $s$ extends to 
$\tilde s \in H^{0}(\f^{-1}(U),\omega^*_X).$ This proves the extendability of $s$ on $U$ and
$\cF$ is reflexive. $W$ being a smooth surface, $\cF$ is locally free.  
$X_0 \la W_0$ being
a conic bundle, there is an embedding $X_0 \hookrightarrow \bP (\F \vert W_0).$ Let $\tilde X$
be the closure in $\bP (\F).$ Then $\tilde X$ is clearly Gorenstein and we claim that $\tilde X$ is a (possibly
singular) conic bundle. 
To see this we let $\pi: \bP (\cF) \la W$ denote the projection and we must prove that there
is no point $w \in W$ such that $\pi^{-1}(w) \subset \tilde X.$ Consider the canonical morphism
$$ \alpha : \f^{*}(\cF) = \f^{*}\f_*(\omega_X^*) \la \omega_X^*$$ 
Let $\cS = {\rm Im}\alpha.$ Then we obtain an embedding
$$ \bP(\cS) \subset \bP(\f^*(\cF)) = \bP(\cF) \times_W X,$$
hence an embedding $\bP(\cS) \subset \bP(\cF).$ It follows that $\tilde X$ is the unique irreducible
component of $\bP (\cS)$ which is mapped onto $W$ by $\pi.$ Assuming the existence of a point
$w \in W$ as above, we have $\bP_2 \simeq \pi^{-1}(w) \subset \bP(\cS).$ If however
$$ p : \bP(\cS) \la W$$
denotes the canonical projection, then, factorising $p$ as $\bP(\cS) \la X \la W,$ it
is clear that $p^{-1}(w)$ cannot be $\bP_2,$ since $\varphi$ us equidimensional [Cu88], contradiction. 
Hence $\tilde X$ is a conic bundle. 
Now there is a birational map $X \rightharpoonup \tilde X,$ which is an isomorphism outside
finitely many curves. Hence $X \simeq \tilde X$ by [Ko89,2.1.13] and $X$ is Gorenstein and
a conic bundle. Note that no component of a fiber of $\f$ is contractible so that [Ko89,2.1.13]
is applicable.    
\sn Now we write $ -4K_W = \Delta + D$ with a nef divisor $D.$ 
\sn (C1) First we consider the case that $e = 0.$ Then $-4K_W \equiv 8C_0.$ Consequently
$\Delta \equiv aC_0$ and $D \equiv bC_0$ with $a + b = 8.$ 
\sn So $\Delta$ consists of $a$ disjoint 
sections. Let $y \in C$ and let $X_y$ be the fiber of $\alpha$ over $y;$ clearly $X_y$ is
reduced. Since $\Delta$ is smooth, every singular conic $\f^{-1}(x), x \in W,$ is a pair of
two different lines. Let $$l = \beta^{-1}(y),$$  $\beta : W \la C$ the projection. Then $l$ meets
$\Delta$ transversally in $a$ points and therefore for $y$ general,  $X_y$ is the blow-up of a Hirzebruch
surface in $a$ points. In particular, $K_{X_y}^2 = 8-a$ for all $y.$ Suppose $X_y$ singular. Consider the projection
$$p: X_y \la l = \bP _1.$$ Since the only singular fibers of $p$ are line pairs, we see that 
$X_y$ has only finitely many singularities. $X_y$ being Gorenstein (because $X$ is Gorenstein),
we conclude that $X_y$ is normal. Let $$ \sigma : \hat X_y \la X_y$$ be the minimal 
desingularisation and
$$\mu : \hat X_y \la \tilde X_y$$ a map to a minimal model. We can arrange things such that 
$\hat X_y \la l$ factors through a map $\tilde X_y \la l$ (just make $\hat X_y \la l$
relatively minimal). We conclude that $\sigma$ contract only parts of fibers of $\hat X_y \la l$
(and
hence $X_y$ has only rational double points as singularities). Since $K_{\hat X_y}^2 = 8-a,$ the birational
map $\mu$ consists of $a$ blow-ups. On the other hand, $X_y \la l$ has exactly $a$ singular
fibers which are line pairs. Therefore $\sigma $ cannot contract any curve, so that $X_y$
is smooth. Hence $\alpha$ is a submersion. In particular $X$ is smooth.
  
\sn (C2) The argument in case $e = -1$ is essentially the same, we thus omit it.

\bn {\bf Remark.} In case $X$ is smooth in the situation of (1.7) and if $-K_X$ is nef, we can 
prove the smoothness of $\alpha$ by direct local calculations, see (4.7).

\bn {\bf 1.8  ÊProposition} \  {\it In (1.7) $-K_X$ is always nef. Moreover the discriminant locus $\Delta$ of the conic bundle
$\varphi$ is - after finite \'etale cover of the base $C$ - of the form $\Delta \equiv \nu C_0,$
where $C_0$ is a section of $W$ with $C_0^2 = 0.$ If $\nu \geq 3$ or with $W = \bP(\cO \oplus L)$
with $L$ a torsion line bundle, then 
 $\varphi$ is analytically a ${\bP_1}-$bundle,
i.e. a conic bundle with discriminant locus $\Delta = \emptyset.$ } 

\bn {\bf Proof.} We make use of the notations of the proof of (1.7). Suppose $\Delta \ne 0.$
We know that $\Delta$ is smooth and that $-K_W$ is nef. 
\sn (1) In a first step we reduce to the case $W = \bP_1 \times C.$ 

(a) If the invariant $e = -1$, take a curve $C_0$ with $C_0^2 = 1.$ By
[At57,lemma 22], $W$ has three \'etale multi-sections $C_i$ of degree 2, which are numerically
equivalent to $2C_0 - F.$ Take one of them, say $C_1$ and perform the base change $C_1 \la A$
to obtain the new ruled surface $W'.$ Then $W'$ has invariant $e = 0.$ 
Hence the case $e = -1$ is reduced to the case $e = 0.$

(b) Since $-K_W$ is nef, (a) implies $e = 0.$ In that case $\Delta \equiv \nu C_0,$ where
$C_0^2 = 0$ and $1 \leq \nu \leq 8.$ In fact, since $-K_F$ is nef, $F$ is a Hirzebruch surface 
blown up in at most 8 points
and therefore 
$$\Delta \equiv \nu C_0 + \mu l,$$
where $l$ is a fiber of $\beta$ (compare the proof of (1.7)). Since on the other hand $-(4K_W + \Delta)$ is nef and
since $-K_W \equiv 2C_0,$ we must have $\mu = 0$ and $1 \leq \nu \leq 8.$ 
We now show that if $\nu /geq 3,$ then we can reduce ourselves to $W = \bP \times C.$ 
If $W = \bP (\cO \oplus L)$ with a topologically trivial line bundle $L,$
then $\Delta$ provides a multi-section, disjoint from the two canonical sections. Hence
$W = \bP_1 \times C$ after a finite \'etale base change. Therefore we may assume that $W$
is a product in that case. If $W = \bP (E)$ with $E$ a nontrivial extension of $\cO$ by
$\cO$, then $\Delta$ provides a multi-section disjoint from the canonical section, so that
after a finite \'etale base change $h: \tilde C \la C,$ the pull-back $\tilde W$ has two
disjoint sections, so that $h^*(E)$ splits. This is impossible.

\bn (2) We consider here the case $W = \bP \times C.$ Let $p_i$ denote the projections of $W$
to $\bP_1$ and $C.$ We consider the fibration $g = p_1 \circ \varphi : X \la \bP_1.$ Its general fiber
$G$ is a$\bP_1-$bundle over an elliptic curve with $-K_G$ nef, hence $G$ has invariant $e = 0$
or $e = -1.$ We can write
$$ \Delta = \bigcup \{x_i\}Ê\times C.$$
Let $C_i = \{x_i\}Ê\times C$ and $G_i = \varphi^{-1}(C_i).$ Then every fiber of $G_i la C_i$
is a reducible conic and thus there exists an unramified 2:1 cover $\tilde C_i \la C_i$ 
such that $\tilde G_i = G_i \times_{C_i}Ê\tilde C_i \la \tilde C_i$ is a $\bP_1-$bundle. 
The map $h: \tilde G_i \la G_i$ is nothing than the normalisation of $G_i.$ By adjunction
$$ K_{G_i}Ê= K_X \vert G_i,$$
hence, $-K_X$ being almost nef, it is clear that $-K_{G_i}$ is nef. If $e_i$ is the invariant
of $\tilde G_i$, it follows as above that $e_i \in \{0,-1\}.$ 
We have the well-known formula (see [Mo82])
$$ K_{\tilde G_i}Ê= h^*(K_{G_i}) - \tilde N, \eqno (*)$$
where $N$ is the non-normal locus (with structure given by the conductor ideal) and
$\tilde N $ the analytic preimage of $N.$ 
Write
$$ h^*(-K_{G_i}) \equiv \alpha C_0 + \beta F,$$
where as usual $C_0$ is a section with $C_0^2 = - e_i$ and $F$ is a ruling line. 
Since $h(F)$ is an irreducible component of a conic in $X,$ it follows
$$ \alpha = h^*(-K_{G_i}) \cdot F = 1.$$
By virtue of $K_{G_i}^2 = K_G^2$ we have
$$ h^*(-K_{G_i})^2 = (C_0 + \beta F)^2 = 2 \beta - e_i = 0,$$ 
in particular $e_i = 0.$ From (*) and 
$$ - K_{\tilde G_i}Ê\equiv 2C_0 + e_i F \equiv 2C_0$$
it follows  $$\tilde N \equiv C_0$$ and
$$ h^*(-K_{G_i}) \equiv C_0.$$ 
Hence $K_{G_i}$ is (numerically) not divisible by 2. Thus $K_G$ is not divisible by 2,
hence $e = e(G) = -1.$ If $C'_0$ and $F'$ are the canonical section resp. a ruling line,
we have $-K_G \equiv 2C'_0 + F'.$ Taking limits yields 
$$ h^*(-K_{G_i}) \equiv 2C_0 + 2F,$$
contradiction. 
\sn Hence $\Delta = \emptyset.$ Now it is clear that every fiber of $\alpha$ is $\bP_1 \times
\bP_1$ and therefore $-K_X$ is nef. 

\bn (3) Next we treat the case $\nu = 1.$ Hence $\Delta$ is a section of $\beta : W \la C$
with $\Delta^2 = 0.$ Then the general fiber of $\alpha : X \la C$ is either 
{\item {(a)} $\bP_1 \times \bP_1$ blown up in one point or 
\item {(b)}Êthe first Hirzebruch surface
$F_1$ blown up in one point, i.e. $\bP_2$ blown up in two points.}
\sn (a) We want to apply (0.5). We start with an irreducible component $B$ of a reducible conic
sitting in a general fiber $F.$ In other words, we consider $\alpha \vert F,$ which is a
$\bP_1-$ bundle over a rational curve blown up in one point and we take a (-1)-curve in a
fiber of $\alpha.$ Since $\varphi$ is a conic bundle, every deformation of $B$ is still a
(-1)-curve in some fiber of $\alpha$ so that we can apply (0.5). We obtain an \'etale cover
$\tilde C \la C$ and a base changed $\tilde \varphi: \tilde X \la \tilde C$ and a birational
morphism $\tau : \tilde X \la X'$ contracting a (-1)-curve in every fiber of $\alpha.$ 
We obtain a submersion $g: X' \la C$ with general fiber $\bP_1 \times \bP_1.$ If we know that
every fiber of $g$ is $\bP_1 \times \bP_1,$ then $-K_{X'}$ is $g-$nef. Since $-K_{X'}$ is
almost nef, we conclude that $-K_{X'}$ is nef, hence $-K_X$ is nef and we are done. This is
certainly the case if $g$ is a contraction of an extremal ray. If $g$ is not an extremal
contraction, we can choose some contraction, say $h: X' \la Z,$ inducing a map $h': Z \la \tilde C.$
It follows that $\dim Z = 2$ and that $h$ is a conic bundle. Since however every $F'$ is a
Hirzebruch surface, it is clear that $h$ must be a $\bP_1-$bundle, and therefore $g$ is
a $bP_1 \times \bP_1-$bundle. 
\sn (b) We proceed in the same way. Now the general fiber of $g$ is $F_1.$ Then either we
can repeat the process by another application of (0.5) or we argue as follows. Since $F' \simeq F_1,$
it is well known that $h$ cannot be an extremal contraction. As in (a) we choose a 
contraction $h: X' \la Z.$ If $\dim Z = 2,$ we conclude as in (a). If $h$ is birationa, then
the general fiber of $h' : Z \la C$ is $\bP_2,$ thus $h'$ is a $\bP_2-$ bundle and $h$ is a
$F_1-$ bundle. Therefore $-K_X$ is nef. 

\bn (4) The case $\nu = 2$ is completely analogous; details are omitted.

\bn {\bf (1.9)} To end the discussion of contractions of fiber type, we must consider the case
${\rm alb}(X) = 2, \dim W = 2.$ Of course we assume $-K_X$ to be almost nef. In that case
$\alpha : X \la {\rm Alb}(X) = A$ has connected fibers [Ka81,Mo87,11.5.3]. Therefore the
map $\beta : W \la A$ has connected fibers, thus it is birational. We claim that $\beta$
is an isomorphism. 

In fact, by (1.6) $-(4K_W + \Delta)$ is nef, $\Delta$ denoting the  discriminant locus of
the "generic" conic bundle $\f$ (cp. (1.5)). Let $h : \hat W \la W$ be
the minimal desingularisation. Since the singularities of $\hat W$ are all rational double points,
we have $K_{\hat W}Ê= h^*(K_W).$ We conclude that $-K_{\hat W}$ is the sum of an effective and a nef
divisor. But $\kappa (\hat W) = 0,$ therefore $\Delta = 0,$ and  $-K_{\hat W} \equiv 0.$ 
So $\hat W$ is a torus, and $\beta \circ h$ and in particular $\beta$ are isomorphisms. 

So $W = A.$ Then (1.6) once again proves $\Delta = 0$ so that $\f$ is analytically a $\bP_1-$bundle 
outside a finite set of $A.$

\bn {\bf 1.10 \ ÊProposition} {\it In the situation of (1.9) $X$ is analytically a $\bP_1-$bundle
over $A.$ In particular $X$ is smooth and $-K_X$ is nef. }

\bn {\bf Proof.} First note that $\f$ is equidimensional, as in (1.7). Let
$$ S = \{a \in A \vert X_a {\rm \  is \ singular}Ê\}.$$
Then $S$ is finite (or empty), (1.5,1.8). So $X \setminus \f^{-1}(S) \la A \setminus S$
is a $\bP_1-$bundle and the same technique as in (1.7) shows that $\f$ is a $\bP_1-$bundle
(noting that no fiber of $\f$ contains a contractible curve since $\f$ is an extremal
contraction). The required torsion freeness of $H^2(A \setminus S, \cO^*)$ follows as in (1.7);
it is equivalent to the torsion freeness of 
$$ H^3(A \setminus S,\bZ).$$  
Now $\varphi$ being a $\bP_1-$bundle, $X$ is smooth and $-K_X$ is nef. 
    
\bn \bn \bn 
\noindent {\medium {2. Birational Contractions}}  
\bn We shall always assume that $X$ is a terminal ${\bf Q}$-factorial threefold with $-K_X$ almost nef.

\bn {\bf 2.1 Proposition}  {\it Let $\f : X \la Y$ be a divisorial contraction. Then
$-K_Y$ is almost nef. }

\bn {\bf Proof.} Let $E \subset X$ be the exceptional prime divisor contracted by $\f.$
If $\dim \f (E) = 0,$ then our claim is obvious; hence we shall assume $\dim \f (E) = 1$
from now on. Let $C = \dim \f (E).$ We only have to show that $K_Y \cdot C \leq 0,$ if $C$
is irrational.
We let $ g: \tilde E \la E$ and $\nu: \tilde C \la C$ be the normalisations and denote $g$
the genus of $\tilde C.$ We obtain a map $p : \tilde E \la \tilde C.$
 \sn Let $h = g \circ \sigma : \hat E \la E.$ Let $f: \hat E \la E_0$ be the minimal
model (note that $\hat E$ is irrational!). Let $C_1 \subset E_0$ be a section with minimal
self-intersection and put $C_1^2 = -e.$ Let $F$ be a general ruling line of $E_0.$ Choose $\lambda$
such that both $\lambda K_X, \lambda K_Y$ are Cartier. Then we have
$$ \lambda K_X = \f^*(\lambda K_Y) + \lambda E,$$
since $\f$ is generically the blow-up of $C.$
It follows that $ h^*(\lambda K_E) = h^*(\lambda K_X \vert E + \lambda E \vert E)$ is Cartier.  
Write $$ h^*(-\lambda K_X) = f^*(\alpha C_1 + \beta F) + \sum a_i A_i,$$
where the $A_i$ are the exceptional components of $f.$
Since $\f$ is generically a blow-up, we see immediately that 
$$ \lambda = \alpha.$$ 
By the same reason we have
$$ h^*(\lambda E \vert E) = f^*(-\lambda C_1 + \gamma F) + \sum b_i A_i.$$
We conclude by adjunction
$$ h^*(\lambda K_E) \equiv f^*(-2\lambda C_1) + (\gamma - \beta) F) + \sum (b_i-a_i)A_i.$$
Now - passing to the level of sheaves - $\omega_{\hat E}$ is a subsheaf of $h^*(\omega_E).$ 
Thus
$$ \omega_{E_0}^{\lambda} = f_*(\omega_{\hat E}^{\lambda}) \subset (f_*h^*(\omega_E^{\lambda}))^{**}.$$
Since $K_{E_0}Ê\equiv -2C_1 + (2-2g-e)F,$ we obtain
$$ -2\lambda C_1 + (\gamma - \beta) F \equiv \lambda K_{E_0}Ê+ \rho F $$
with $\rho \geq 0.$ Squaring yields
$$ -4\lambda^2e + 4\lambda \beta - 4\lambda \gamma = 8 \lambda^2 (1 - g) - \lambda g,$$
hence $- \lambda e + \lambda \beta - \lambda \gamma \leq 0.$ This implies $\beta + \gamma \geq 0.$

\bn {\bf 2.2 Proposition} {\it Let $\f : X \rightharpoonup X^+ $ be a flip. If $-K_X$ is almost
nef, then so does $-K_{X^+}.$} 

\bn  {\bf Proof.} Let $E \subset X$ and $E^+ \subset X^+$ be the exceptional sets so that
$\f : X \setminus E \la X^+ \setminus E^+$ is an isomorphism. Both $E$ and $E^+$ consist of
finitely many rational curves, so we do not have to care about curves in $E^+.$ 
Therefore it is sufficient to show the following:

if $C \subset X$ is an irreducible curve, $C \not \subset E,$ and if $C^+ \subset X^+$ denotes
its strict transform, then $K_{X^+}Ê\cdot C^+ \leq K_X \cdot C.$ 

Choose a desingularisation $g: \hat X \la X$ such that the induced rational map $h : \hat X \la X^+$
is a morphism. Then one has
$$ K_{\hat X}Ê= g^*K_X + \sum \lambda_i E_i $$
and
$$ K_{\hat X}Ê= h^*K_{X^+}Ê+ \sum \mu_i E_i,$$
where the $E_i$ are the exceptional components of $g.$ 
Then by [KMM87,5.1.11] we have $\lambda_i \leq \mu_i$ from which our inequality is clear.

\bn {\bf 2.3 Proposition} \ ÊÊ{\it Let $X$ be a smooth projective threefold with $-K_X$ nef and
positive irregularity $q(X).$ Let $\f: X \la Y$ be a divisorial contraction. Then $\f$ is the
blow-up of a smooth curve $C \subset Y$ and $-K_Y$ is almost nef. If $-K_Y$ is not nef, then
$C \simeq \bP_1$ with normal bundle $N_{C \vert Y}Ê\simeq \cO(-2) \oplus \cO(-2).$} 

\bn {\bf Proof.} [DPS93]
\bn The exception described in (2.3) is the reason why  we introduce the notion "almost nef". In the
end it will turn out that this exception does not happen. If $q(X) = 0$ then the exception 
might very well occur.

\bn \bn \bn 
{\medium {3. The Main Theorem}}

\bn Here we begin studying backwards: we start with a smooth object with $-K$ nef and ask how
we can modify without destroying this property.
\bn {\bf 3.1  Proposition}  {\it Let $Y$ be a smooth projective threefold with $-K_Y$ nef. Let
$\beta : Y \la A$ be its Albanese map to the abelian surface $A.$ We assume that $\beta$ is
a submersion. Let $X$ be a terminal threefold and let $\f: X \la Y$ be a
divisorial contraction. Then $-K_X$ is not almost nef.}

\bn {\bf Proof.}  We may assume that $\kappa (Y) = - \infty,$ otherwise our claim is
obvious. By our assumption $\beta$ is a $\bP_1-$bundle analytically. Let $E$ be the exceptional
divisor of $\f.$ 
\sn (a) First let $\dim \f (E) = 0.$ We can write
$K_X = \f ^*(K_Y) + \mu E$ for some positive rational number $\mu.$ Notice that $X$ might not be smooth, even not
Gorenstein as examples using weighted blow-ups, say in $\bP_3$, show. 
First of all we have $$ K_X^3 = K_Y^3 + \mu ^3 E^3. $$
Since $K_Y^3 = 0 $ and $E^3 > 0,$ we conclude $K_X^3 > 0,$ hence $-K_X$ cannot be nef.
If $-K_X$ is almost nef, there is a rational curve $C \subset X$ with $K_X \cdot C > 0.$
Then $\varphi (C)$ must be a fiber of $\beta,$ namely the fiber containing $p = \varphi(E).$
In particular $C$ is the unique curve in $X$ with $K_X \cdot C > 0.$ 
Observe that after a possible \'etale base change, we may assume $Y = \bP (V)$
with a rank 2-bundle $V$ on $A.$ Since $-K_Y$ is nef, $V$ is  numerically flat (after another base change) [CP91,DPS94] and
thus we have an exact sequence 
$$0 \la L_1 \la V \la L_2 \la 0$$
with flat line bundles $L_i.$ In particular $K_Y^2 = 0.$ 
Now take a general smooth surface $S$ through $p$ and let $\hat S$ be its strict transform in $X.$ Then
$$ K_X^2 \cdot \hat S = K_Y^2 \cdot S + 2\varphi^*(K_Y) \cdot E \cdot \hat S + \mu^2 E^2 \cdot \hat S = \mu^2 E^2 \cdot \hat S < 0.$$
Hence $-K_X \vert \hat S$ cannot be nef; on the other hand $\hat S$ does not contain $C,$
contradiction.   
\sn (b) If $\dim \f (E) = 1,$ choose a general curve $B \subset A.$ Let $\hat B$ its
preimage under $\beta \circ \f.$ Then $-K_{\hat B \vert B}$Ê is almost nef, hence nef and
therefore $\hat B \la B$ is a submersion (4.4). This proves $\beta \f (E) = 0.$
If $-K_X$ is almost nef, then an argument as in (a) shows that $-K_X$ is nef. But in that
case simple numerical calculations give a contradiction, see (4.19) for the details in a
slightly more general situation.

\bn {\bf 3.2 Proposition}  {\it Let $Y$ be a smooth projective threefold with $-K_Y$ nef, $\beta : Y \la A$
the Albanese to the elliptic curve $A.$ Suppose
that $\beta$ is a submersion. Let $\f : X \la Y$ be the blow-up of a point or a smooth curve $C.$ If 
$-K_X$ is nef, then $\f$ cannot be the blow-up of a point. If $\f$ is the blow-up of $C,$ then 
$C$ is an \'etale multi-section of $\beta $ 
and $\alpha$ is smooth.}

\bn {\bf Proof.} The first claim, that $\f$ is not the blow-up of a point, is obvious since we have
$K_X^3 = K_Y^3 = 0.$ So assume that $\f$ is the blow-up of the curve $C.$ If $ \dim \beta (C) = 1,$ then our claim 
follows from the more general proposition (4.11), therefore we shall assume $\dim \beta (C) = 0,$ so that
$C $ is contained in a fiber $F$ of $\beta.$ The case $K_F^2 > 0$ is treated in (4.11), too.
Hence it remains to consider the case $K_F^2 = 0$. We may assume $\kappa (X) = - \infty.$ Then either 
\sn (1) $F$ is a $\bP_1-$bundle over an elliptic curve with invariant $e \leq 0$ or
\sn (2) $F$ is $\bP_2$ blown up in nine sufficiently general points. 

\bn (1) In this case $\beta$ factors by (0.4) in the following way
$$ Y \buildrel \gamma \over {\la} Z \buildrel \delta \over {\la} A $$
with $\gamma$ a $\bP_1-$bundle and $\delta$ an elliptic bundle. Hence $Z$ is hyperelliptic. 
Then we perform an \'etale base change $\tilde Z \la Z$ with $\tilde Z$ an abelian surface  and conclude
easily by applying (3.1). 

\bn (2) Here we have a factorisation 
$$ Y \buildrel {\gamma_1} \over {\la} Y_1 \la \ldots \buildrel {\gamma_k} \over {\la} Y_k \buildrel {\delta}Ê\over {\la} Z 
\buildrel {\epsilon} \over {\la} A $$
or
$$ Y \buildrel {\gamma_1} \over {\la} Y_1 \la \ldots \buildrel {\gamma_k} \over {\la} Y_k 
\buildrel {\rho} \over {\la} A$$    
with $\gamma_j$ blow-ups of \'etale multi-sections, $\delta$ and $\epsilon$ both $\bP_1-$bundles and $\rho$ a $\bP_2-$bundle. 
Let $C_j$ be the image of $C$ in $Y_j.$ Let $E_j \subset Y_{j-1}$ be the exceptional divisor 
of $\gamma_j$ and let
$B_j$ be the center of $\gamma_j$ so that $E_j = \gamma_j^{-1}(B_j).$ 
Now by the computations of [DPS93,p.234-235] and of the proof of (4.11) below we have, in
the notations of (4.11) that $b + \mu = 0,$ which is to say that $K_Y \cdot C = 0.$ Hence
$K_F \cdot C = 0$ and $C$ is an elliptic curve. Moreover
there is an index $j$ such that $C_{j-1}Ê\cap E_j \ne \emptyset,$ i.e. $C_j \cap B_j \ne 0.$ 
\sn Every $B_j$ is an elliptic curve; we now check that the normal bundle $N_{B_j \subset Y_j}$ is flat (hence the
ruled surface $E_j$ has invariant $e = 0).$ In fact, we see inductively that $-K_{Y_{j-1}}$ is nef and
that $K_{Y_{j-1}}^3 = 0.$ 
As in [DPS93,p.234] and (4.11) we write
$$ -K_{Y_{j-1}} \vert E_j \equiv C_0  + bf,$$ and
$$ N_{E_j \vert Y_{j-1}} \equiv -C_0  + \mu f,$$
where $f$ is a ruling line of $E.$ Then we have by [DPS93]:
$$ b+\mu = 2b-e.$$ 
On the other hand we have by the proof of (4.11) that
$$b+\mu = 0 \ {\rm and }Ê\ b = {e \over 2}.$$ 
Since $e \geq -1,$ we conclude $e \geq 0,$ hence $b = \mu = e = 0.$ So $N_{B_j \vert Y_j}$ is flat.

\sn With this last observation our claim now clearly follows from the

\bn {\bf 3.3.a Sublemma} {\it Let $Z$ be a smooth projective threefold,  $B \subset Z$ be a smooth elliptic curve with flat
normal bundle $N_B$ and $\psi : Y \la Z$ be the blow-up of $B.$ Denote $E = \psi^{-1}(B)$ the exceptional divisor of $\psi.$ 
Let $C \subset Y$ be a smooth curve and $\varphi : X \la Y$ be the blow-up of $C.$ Assume $-K_Y$ nef
and $K_Y^3 = 0.$ Then $-K_X$ is not nef unless $C \cap E = \emptyset.$ }

\bn {\bf Proof.} We use analogous notations as in part (2) of the proof of (3.2) and have 
by our assumptions (cp. [DPS93,p.234-235], (4.11)):
$$ -K_Y \vert E \equiv 0; \ K_E \equiv -2C_0.$$
Let $\hat E$ be the strict transform of $E$ in $X.$ Then 
$$ K_X \vert  \hat E = \varphi^*(K_Y \vert E) + D,$$
where $D$ is an effective divisor supported exactly on the exceptional set of $\hat E \la E.$
Now suppose that $C \cap E \ne \emptyset.$ Then $- K_X \vert \hat E \equiv \varphi^*(C_0) - D$,
and, $D$ being non-zero, we conclude 
$$ K_X^2 \cdot \hat E = D^2 < 0,$$
so that $-K_X$ cannot be nef. 

\sn This finishes the proof of both the Sublemma and (3.2). 
\sn In the proof of (3.3) we will see that (3.2) remains true also if we only suppose 
$-K_X$ to be almost nef, but this turns out to be much more complicated.

\bn We are now in the position to prove the main result of this paper.

\bn {\bf 3.3Ê Theorem}  {\it Let $X$ be a smooth projective 3-fold with $-K_X$ nef. Then
the Albanese map $\alpha : X \la A$ is a surjective submersion.} 

\bn {\bf Proof.} We know already by [DPS93] that $\alpha$ is surjective. Of course we may
assume that $q(X) > 0.$ If $K_X$ is nef, then by (1.2) $K_X \equiv 0$ and it is well-known
(see e.g. [Be83]) that after a finite \'etale cover $X$ is a product of a torus and a K3-surface
or $X$ is a torus. Then our assertion is clear. So we shall assume that $K_X $ is not nef.
Then $\kappa (X) = - \infty,$ hence $X$ is uniruled, $q(X) \leq 2$ and there exists an 
extremal contraction $$\f: X \la X_1.$$ We have a factorisation $\alpha = \beta \circ \f$
with $\beta : X_1 \la A$ the Albanese of $X_1$ (of course $X_1$ might be singular).  
\sn (1) First assume that $\dim X_1 < \dim X.$ Then we conclude by (1.3),(1.7) and (1.10).

\sn (2) Now suppose that $\dim X_1 = \dim X$ {\rm and} that $-K_{X_1}$ is nef. Let $E$ be
the exceptional divisor of $\f.$ Then $\dim \f(E) = 1;$ [DPS93,3.3]; otherwise $(-K_{X_1})^3 > 0$ so that
$-K_{X_1}$ would be big and nef, hence $q(X_1) = 0$ by [KoMiMo92]. Hence $X_1$ is smooth and
by induction on $\rho (X)$ we conclude that $\beta$ is a submersion. Then $\alpha $ is smooth
by (3.1) and (3.2).

\sn (3) Finally we deal with the case that $\dim X_1 = \dim X$ and that $-K_{X_1}$ is not nef.
By [DPS93] this happens exactly when the exceptional divisor $E$ is mapped to a smooth rational
curve $C \subset X_1$ with normal bundle $N_C = \cO(-2) \oplus \cO(-2).$ Moreover $K_X^2 \cdot F = 0$
for every fiber $F$ of $\alpha.$  
We must show that this special situation cannot occur.

\noindent To this extend we perform Mori's minimal model programme and obtain a sequence
$$ X \la X_1 \la X_2 \la ... \la X_k \la X_{k+1}$$
of extremal contractions $\f_i : X_i \la X_{i+1}$ resp. flips $\f_i : X_i \rightharpoonup X_{i+1}$
such that $$\dim X_k = 3,  \dim X_{k+1}Ê\leq 2.$$
In order to simplify notations we let $Y = X_k$ and $Z = X_{k+1}.$ Furthermore let
$f = \f_k.$ The map $\alpha$ clearly induces maps $\beta : Y \la A$ and $\gamma : Z \la A$
such that $\beta = \gamma \circ f.$ 

\noindent By (2.1) and (2.2) $-K_Y$ is almost nef. Hence by (1.3),(1.7), (1.8) and (1.10) $Y$ is smooth,
$-K_Y$ is nef and $\beta$ is a submersion. 
 It follows that $\f_{k-1} : X_{k-1}Ê\la Y$ cannot be a flip, so it has
to be a divisorial contraction. If $\dim A = 2,$ we apply (3.1) to conclude that
$-K_{X_{k-1}}$ cannot be almost nef which contradicts the nefness of $-K_X$ via (2.1) and
(2.2).

\noindent Therefore we are left with the case that $\dim A = 1.$ Then either 
 \bn {\bf Case I} \ Ê{\rm $\beta$ is the contraction of an extremal ray, in particular $Z = A,$} or

\bn {\bf Case II} \ ÊÊ{\rm $\dim Z = 2$}.
 
\bn We are going to show that in both cases the sequence
$$ X \la ... \la X_{k-1}Ê\la Y$$
consists of blow-ups of \'etale multi-sections over
$A.$ To prove this, we proceed step by step starting with $X_{k+1} = Y$ and we are allowed to perform
\'etale base changes on $A.$

\bn {\bf Case I}
\smallskip \noindent By (1.3), $\beta$ is a submersion so that $\beta$ is a
$\bP_2-$ or $\bP_1 \times \bP_1-$ bundle (0.4). In the second subcase we can reduce by a base
change to Case II, applying [CP91,7.2]Ê($\beta$ becomes a $\bP_1-$bundle over
a $\bP_1-$bundle). For simplicity of notations let
$W = X_{k-1}$.
We can write $$ Y = \bP (E_0)$$ 
with a 3-bundle $E_0$ over $A.$  The nefness of $-K_Y$ is equivalent
to saying that $$ E_0 \otimes {1 \over 3} {\rm det}ÊE^*_0 $$
is nef or that $E_0$ is semi-stable. By another base change and normalisation, taking into 
account [At57], we have the following situation. There are exact sequences
$$ 0 \la \cO \la E_0 \la F_0 \la 0 \eqno (S_1)$$
and
$$ 0Ê\la L_1 \la F_0 \la L_2 \la 0 \eqno (S_2) $$
with a 2-bundle $F_0$ and topologically trivial bundles $L_i$ on $A.$
Therefore we have a distinguished surface $$ \bP := \bP(F_0) \subset Y.$$
Of course the sequences $(S_i)$ might not be unique. 

\bn (A)  
We are now going to investigate the structure of $g: W \la Y$ and will show that $g$ is the blow-up
of a "canonical" section coming from some sequence $(S_i).$ First note that $g$ cannot be small 
since $-K_Y$ is nef. So
let $E$ be the exceptional divisor of $g.$
We claim that 
$$ \dim g(E) = 1.$$
Suppose $\dim g(E) = 0.$ Then we argue similarly as in the beginning of (3.1). Namely, if $-K_W$ is nef, then $K_W^3 = 0,$ hence $K_W = g^*(K_Y) + \mu E$ 
with positive $\mu,$ easily (as before) gives
$$ E^3  = 0,$$
which is absurd. Hence $-K_W$ is not nef but almost nef. Since $g$ is a weighted blow-up (of type $(1,a,b)$ 
with relatively prime positive integers $a$ and $b$, i.e. the blow-up of the ideal $(x,y^{a},z^b) $ in suitable coordinates, cp.Ê[Ka96]), it is immediately calculated that
$-K_W$ is still relatively nef over $A.$  Since $-K_W$ is not nef, we find an irreducible curve $C$ 
such that $K_W \cdot C > 0.$ Then $C$ maps onto $A$,
so that $C$ is irrational. Hence $-K_W$ is not almost nef. 
\noindent So $g(E) $ is a curve $D.$ We are going to show that $D \la A$ is \'etale so that
$D$ is a smooth elliptic curve, $W$ is smooth and $g$ the ordinary blow-up. Of course $g$ is
generically the blow-up of the smooth curve $D.$

\bn We will distinguish three different cases according to the position of $D$ and $\bP.$

\sn (a) $D \cap \bP$ is a finite non-empty set.

\noindent 
Let $\hat \bP$ be the strict transform of $\bP$ in $W.$ By abuse of notation we will not
distinguish between $g$ and $g \vert E.$ Let $C_0 \subset \bP$ be a curve with $C_0^2 = 0$
such that $C_0$ and a ruling line $F$ generate the cone of curves. Let $E' = E \cap \hat \bP.$
Then
$$ - K_W \vert \hat \bP = g^*(-K_Y \vert \bP) - E' = g^*(-K_{\bP} + N_{\bP})Ê- E' \equiv
g^*(3C_0) - E';$$
here $N$ denotes the normal bundle.
If $\hat \bP$ happens to be singular, we pass to a desingularisation, so that we shall assume now
$\hat \bP$ to be smooth. It is actually sufficient to consider the case where $\hat \bP \la
\bP$ is the blow-up of one simple point; the other cases will factorise over this case.
We know that $g^*(3C_0) - E'$ must be almost nef. On the other hand
we have
$$ (g^*(3C_0) - E')^2 = -1 ,$$
hence $g^*(3C_0) - E'$ is not nef. But clearly $g^*(3C_0) - E'$ is nef on every rational curve 
of $W,$  contradiction. 

\sn (b) $D \subset \bP.$

\noindent In this case $D$ is locally a complete intersection curve, so that we know a
priori that $g$ is the blow-up of $D$ (0.6). We therefore know $\hat \bP \simeq \bP.$
If $D\subset \bP$ is a ruling fiber, we see immediately that $-K_W \vert \hat \bP$ is not almost nef,
so assume that $D$ is a multi-section of $\bP.$ Since $-K_W \vert \hat \bP \equiv g^*(3C_0) - D,$
we conclude, identifying $\hat \bP$ and $\bP$ and writing $D \equiv aC_0 + bF$, that
$$ 3C_0 - D \equiv (3-a) C_0 - bF$$
must be almost nef. By virtue of
$$ (3C_0 - D) \cdot F \geq 0$$ and
$$ (3C_0 - D) \cdot C_0 \geq 0,$$
we deduce that $D \equiv aC_0$  with $a \leq 3.$ After another base change $D$ becomes
a section and must be of the form $D = \bP(L_i),$ using the sequence $(S_2)$ (if $F_0$ splits,
then of course we need a suitable choice of $L_i$).

\sn (c) $D \cap \bP = \emptyset.$

\noindent Now $D$ is a multi-section of $\beta $, let $h : D \la A$ denote the restriction of
$\beta.$ Then $D$ provides a section of $\bP(h^*(E_0)),$ disjoint from $h^*(\bP).$ Thus 
$h^*E_0 = \cO \oplus h^*(F_0).$ Let $\zeta \in H^1(A,F_0^*)$ denote the 
extension class defining $(S_1).$ By the above splitting it follows $h^*(\zeta) = 0.$
On the other hand the restriction map
$$ h^* : H^1(A,F_0^*) \la H^1(D,h^*(F_0^*)) $$
is injective since $\cO_A $ is a direct summand of $h_*(\cO_D)$ (we may assume that $D$ is
smooth). Therefore sequence $(S_1)$ already splits and we conclude that $D = \bP(\cO).$

\bn (B) Now we have completely determined the structure of $g;$ it is (after base change) the
blow-up of one of the canonical section of $Y$ coming from $(S_1)$ or $(S_2).$ Note also that
clearly $-K_W$ is nef. Now we proceed with the next contraction
$\f_{k-2}: X_{k-2}Ê\la X_{k-1},$ which we rename $g_1: W_1 \la W.$ We proceed in the same
way as before, the arguments being similar. Note that $\bP$ survives (as strict transform)
in $W_1$; we denote the transform again by $\bP.$  First we show
$$ \dim g_1(E') = 0$$ as before 
where $E' \subset W_1$ again denotes the exceptional divisor or in case of a blow-up of the {\it smooth} point $p$ 
we have the following geometric argument. Let $\hat F$ be the fiber component of $\beta \circ g \circ g_1$ such that $p \in
g_1(\hat F).$ 
Then $\hat F \simeq \bP_2(x,p),$ the blow-up of $\bP_2$ at $x$ and $p.$ First we shall assume
that $x$ and $p$ are not infinitesimally near. Then we choose an irreducible cubic
$C \subset \bP_2$ passing through $x$ and $p$ and having multiplicity 2 at $p.$ Let $\hat C$
be its strict transform in $\hat F.$ Then - with $A = \hat F \cap E$ -  
$$\hat C \in \vert -K_{\hat F} - A \vert.$$
Noticing that $$-K_{W_1}Ê\vert \hat F = - K_{\hat F} - A,$$
we conclude, using $\hat C,$ that $-K_{W_1}Ê\vert \hat F$ is nef, and therefore $-K_{W_1}$ is relatively nef
with respect to $\beta \circ g \circ g_1.$ In particular
$$ - K_{W_1} \cdot C \geq 0$$
for all rational curves $C \subset W_1.$ Since $-K_{W_1}$ is almost nef, it is actually nef.
But $K_{W_1}^3 = 1,$ since $K_W^3 = 0,$ contradiction. 

\noindent If $p$ is infinitesimally near to $x$, then $-K_{\hat F} - A$ is no longer nef,
so we argue as follows. We consider the $\bP_1-$bundle $E \la D$ and let $\hat E \subset W_1$
be its strict transform. The normal bundle $N_{D \vert W}$ is a flat vector bundle. Hence
$K_E \equiv -2C_0$ and $N^*_{E \vert W}Ê\equiv C_0.$ Thus
$$ -K_W \vert E \equiv C_0.$$
Now $$ -K_{W_1}\vert \hat E = g_1^*(-K_W) \vert \hat E - 2E' \vert \hat E \equiv g_1^*(C_0) -
2l $$
where $l = E' \cap \hat E.$ Here we have used $p \in E,$ which follows from the fact that $p$
and $x$ are infinitesimally near. We conclude that $g_1^*(C_0) - 2l$ is almost nef.
Now take a section $C \in \vert C_0 \vert $ resp. $C \in \vert C_0 + F \vert $ with a ruling
line $F.$ Then, denoting $\hat C$ the strict transform in $\hat E,$ 
$$ g_1^*(C_0) - 2l \cdot \hat C < 0,$$
contradiction.    

\bn Therefore we know that $g_1$ is not the blow-up of a point, hence it must be centered 
at a curve $D_1.$ First suppose $D_1 \subset E.$ Then $g_1$ is the blow-up of
$D_1$ (0.6). We identify $\hat E $ with $E.$ Then we have
$$ -K_{W_1}Ê\vert E = -K_W \vert E - D_1 \equiv C_0 - D_1.$$
So $C_0 - D_1$ is almost nef. We conclude easily that $D_1 \equiv C_0.$ So $D_1$ is a section
of $W_1 \la C.$ If $D_1 \not \subset E,$ we consider $D \cap \bP$ and conclude as in (A),
distinguishing the cases $D \cap \bP$ finite, empty or $D \subset \bP.$ Again $-K_{W_1}$ is
nef.  

\noindent In the next step we have to consider $g_2: W_2 \la W_1$ and again have to rule out the
blow-up of a point. Here $\hat F = \bP_2(x_1,x_2,p)$ and it is convenient in the case of
general position to choose a line $l_i
\subset \bP_2$ such that $x_1,x_2 \in l_1$ and $p \in l_1 \cap l_2.$ Then 
$$ \hat l_1 + \hat l_2 + \hat l_3  \in \vert -K_{\hat F} - E \vert ,$$
from which the nefness of $-K_{W_2}Ê\vert \hat F$ is an immediate consequence. In the infinitesimal
near case we argue as before. 

\bn (D) Continuing this way we can handle 5 steps (afterwards the linear system $\vert -K_{\hat F}
- E \vert = \emptyset.)$ In every fiber $F_4 $ of $W_4 \la A$ at most two points can be
infinitesimally near, otherwise $-K_{F_5}$ would no longer be nef. Mapping all the centers
of the blow-ups $\f_i$ to $Y$, we therefore obtain at least 4 disjoint multi-sections of
$ Y = \bP(E_0) \la A.$ Comparing with $(S_i),$ we conclude that for a suitable choice of
$(S_i)$ and after possibly substituting $E$ by $E \otimes L,$ with $L$ topologically trivial,
we have either
$$ F_0 = \cO \oplus \cO$$ 
with all the sections to be blown up in $\bP(F_0),$ or 
$$ E_0 = \cO^{\oplus 3}.$$ 
But the first case cannot occur: looking again at a fiber $F_4,$ we then would find 4 points
in $\bP_2$ (to be blown up) on a line which is not possible since $-K_{F_4}$ is nef.
So we have $Y = \bP_2 \times A.$

\bn (E) Now our claim follows very easily: assume that we have done already $j$ steps,
i.e. we have blown up only sections of the form $x_i \times A.$ Then the result $W_{j-1}$
is of the form $$W_{j-1}Ê= \bP_2(x_1,...,x_j) \times A.$$Ê
Now suppose that $g: W_j \la W_{j-1}$ is the blow up of a point $p = (p_1,p_2) \in
W_{j-1}Ê\times A.$ Then let $B = p_1 \times A.$ We see immediately that
$$K_{W_j}Ê\cdot B > 0,$$
contradicting the almost nefness of $-K_{W_j}.$ 
So $g$ contracts a divisor to a curve $C.$ Choose a generic smooth point $p \in C$ and
define $B$ as before. If $C \ne B,$ the same computation as above yields a contradiction,
hence $C$ is as claimed. 

\bn (F) Conclusively $X \la Y$ is the blow-up of \'etale (multi-)sections so that 
$\f : X \la X_1$ cannot be the blow-up of a rational curve. This finishes Case I.

\bn {\bf Case II} 
\sn This case is done partly in the same way, partly reduced to Case I. Note that $-K_Z$ is nef, hence it is either
a hyperelliptic surface, in which case we pass to an abelian 2-sheeted cover of $Z$ so that
we reduce to the case $\dim A = 2.$ Or [CP91] $Z$ is
a $\bP_1-$bundle over $C$, moreover $Z = \bP(E)$ with a semi-stable rank 2- bundle $E$
on $A.$ After passing to a $2:1-$cover of $A,$ the bundle $E$ is flat. If $\psi : Y \la Z$
is a $\bP_1-$bundle, it is given by $Y = \bP(V)$ with a 2-bundle $V$ on $Z$ and it is
clear that $-K_Y$ is nef since it is almost nef. Then we can proceed
in the same way as in Case I. So suppose that $\psi $ is a proper conic bundle.
By (1.8) $-K_Y$ is nef. Note that $Y \la A$ is a submersion since $-K_Y$ is nef.  
Then, using (0.4) we perform another base change
to reduce our situation to Case I.

\sn This finishes the proof of the Main Theorem. 

\bn The proof of the main theorem actually gives a more explicit description of the Albanese map.

\bn {\bf 3.4 Corollary}  {\it Let $X$ be a smooth projective threefold with $-K_X$ nef. 
Let $\alpha: X \la A$ be the Albanese. 
{\item {(1)} If $\dim A = 2,$ then $X$ is a $\bP_1-$ bundle over $A.$
\item {(2)} If $\dim A = 1,$ then there exists a sequence of blow-ups $\f_i: X_i \la X_{i+1}, 
0 \leq i \leq r,$ with $X_0 = X$ and inducing maps $\alpha_i : X_i \la A$  such that
{\itemitem {(a)}all $X_i$ are smooth, all $-K_{X_i}$ are nef, $\f_i$ is the blow up of a smooth curve $C_i$ and
$C_i$ is an etale multi-section of $\alpha_{i+1};$
\itemitem {(b)} the induced map $\alpha_{r+1}Ê: X_{r+1} \la A$ is a $\bP_2-$ bundle or a
$\bP_1 \times \bP_1-$ bundle or it factors as $h \circ g$ with $g : X_{r+1}Ê\la Y $ a
conic bundle and $h: Y \la A$ is a $\bP_1-$ bundle.}}}  

\bn {\bf 3.5 Corollary}  {\it Let $X$ be a smooth projective threefold with $-K_X$ nef. 
Let $\alpha : X \la A$ be the Albanese and assume $\dim X = 1.$ Then there exists a finite
etale cover $\tilde X \la X$ induced by a finite etale cover $\tilde A \la A$ such that the
following holds. There exists a finite sequence of blow-ups of sections over $\tilde A,$
say $\tilde X \la \tilde X_1 \la ... \la \tilde X_{r+1}$ such that the induced map
$\alpha_{r+1}Ê: \tilde X_{r+1}Ê\la \tilde A$ is $\bP_2-$bundle or a $\bP_1 \times \bP_1-$
bundle over $\tilde A.$ In the first case $\tilde X_{r+1}Ê= \bP(E)$ with a semi-stable vector
bundle of rank 3 on $\tilde A.$ In the second case $\alpha_{r+1}$ is the contraction of an
extremal ray (hence $\rho(\tilde X_{r+1}) = 2$) or $\alpha_{r+1}$ factorises as $ \alpha_{r+1} =
g \circ f,$ where $f: \tilde X_{r+1}Ê\la S$ is a $\bP_1-$bundle and $S = \bP(F)$ with 
$F$ a semi-stable rank 2 - bundle over $\tilde A.$}

\bigskip \bigskip \bigskip
{\medium 4. The relative case}

\bn In this section we want to consider the following situation. Let $X$ be a smooth projective
manifold of dimension $n$ and $Y$ a projective manifold, $\dim Y \geq 1.$ Let $\f : X \la Y$
be a surjective map. Assume that $$ -K_{X \vert Y}Ê= \omega^{-1}_{X \vert Y}Ê= \omega^{-1}_X 
\otimes \f^*(\omega_Y)$$
is nef. What can one say about the structure of $\f?$ Our previous situation of the last three
sections is the special case when $\dim X = 3,$ $Y$ is abelian and $\f$ the Albanese. We shall fix the above
notations for the entire section. Miyaoka has shown in [Mi93] that $\XY$ is never ample. His
proof works for all ground fields, even not algebraically closed. For algebraically closed
fields of characteristic 0 the statement can be improved, the proof being much easier:

\bn {\bf 4.1 Proposition} {\it Suppose $\omega_{X \vert Y}^{-1}$ nef. $\XY$ is not big, i.e. $(\XY)^n = 0.$ }

\bn {\bf Proof.} We proceed by induction on $d = \dim Y.$ Suppose that $\XY$ is nef and big. 
By Kawamata-Viehweg vanishing we obtain:
$$ 0 = H^1(X,\XY \otimes \omega_X) = H^1(X,\f^*(\omega_Y)).$$ 
The Leray spectral sequence yields $$H^1(Y,\omega_Y) = 0,$$
which gives a contradiction in case $d = 1.$ 
\sn Now suppose that the claim holds for values of
$\dim Y$ smaller than $d.$ Take a smooth very ample divisor $Z \subset Y$ such that $W :=
\f^{-1}(Z)$ is smooth. Then $\XY \vert W = \omega^{-1}_{W \vert Z}$ is nef, hence by
induction $\omega^{-1}_{W \vert Z}$ is not big. Hence
$$ 0 = (\omega^{-1}_{W \vert Z})^{n-1} = (\XY)^{n-1}Ê\cdot W.$$
But clearly $\kappa(\XY \vert W) = n-1$ for general choice of $Z.$ This is a contradiction. 

\bn As the referee pointes out, (4.19) remains true if $-(K_X + \Delta)$ is nef as long as 
$(X,\Delta)$ is log terminal in the sense of Kawamata [KMM87].

\bn {\bf 4.2 Proposition}  {\it  Suppose that $\XY$
is nef and that the general fiber has Kodaira dimension $\kappa (F) \geq 0.$ Then 
$\kappa (F) = 0, \f$ is smooth and locally trivial and $\XY$ is a torsion line
bundle. }

\bn {\bf Proof.} Since $-K_F = \XY \vert F$ is nef and $\kappa (F) \geq 0,$ it follows that
$K_F$ is torsion. Choose a positive integer $d$ such that $dK_F = \cO_F.$ Then 
$\f_*(\omega_{X \vert Y}^{\otimes d})$ is of rank 1. Viehweg has shown that $\f_*(\XXY^{\otimes d})$
is weakly positive, see e.g. [Mo87] for definition and references. Since the natural injective
map
$$ \f_*(\XXY^{\otimes d})Ê\la (\f_*(\XY^{\otimes d}))^{**}$$
is generically surjective, it turns out that $$(\f_*(\XXY^{\otimes d}))^{**}$$ is weakly positive,
too [Mo87,5.1.1(b)]. But this last sheaf is invertible, and for invertible sheaves the notions
of weak positivity and pseudoeffectivity are equivalent [Mo87,p.293]. Therefore 
$(\f_*(\omega_{X \vert Y}))^{**}$ is pseudoeffective, i.e. numerically equivalent to a limit of
effective $\bQ-$divisors, and so does its pull-back to $X.$ Via the generically surjective map
$$ \f^*(\f_*(\XXY^{\otimes d}))^{**}Ê\la \XXY^{\otimes d},$$
we conclude that $\XXY^{\otimes d}$  is pseudoeffective. We claim that 
$$ \XXY \equiv 0.$$
In fact, take an ample divisor $H$ on $X.$ By pseudoeffectivity we have $\XXY \cdot H^{n-1}Ê\geq 0$
while by our nefness assumption we get the reversed inequality. Hence $\XXY \cdot H^{n-1}Ê= 0$
which easily implies our claim. By Corollary 1.2 in [Ka85] we deduce
$$ \kappa(\XXY) \geq \kappa (F) = 0.$$ 
Hence $\XXY$ is torsion. Finally Theorem 4.8 in [Fu78] shows that $\f$ is smooth and 
locally trivial. 

\bn {\bf Remark.} The hypothesis $\kappa (F) \geq 0$ in (4.2) can be (formally) weakened to
$K_F \cdot H^d \geq 0$ for some ample divisor $H$ on $F, \dim F = d.$ In that case, keeping in 
mind that $-K_F$ is nef, we get $K_F \cdot H^d = 0,$ which implies $K_F \equiv 0.$
Then Kawamata's result [Ka85,8.2] shows that $\kappa (F) = 0$ so that $K_F$ is torsion. 

\bn {\bf 4.3 Proposition}  {\it Let $X$ be a terminal threefold and $\f : X \la Y$ a surjective
morphism to a normal projective $\bQ-$Gorenstein surface. Assume that $-K_{X \vert Y}$ is
nef. Then $\dim \{y \in Y \vert X_y {\rm \ is \ singular}Ê\}Ê\leq 0.$ } 

\bn {\bf Proof.} Let $C \subset Y$ be a general irrational hyperplane section. Then $X_C := \f ^{-1}(C)$
is smooth and $-K_{X_C \vert C}$ is nef. Now apply the following proposition (4.4).

\bn {\bf 4.4 Proposition} {\it Let $f : S \la C$ be a surjective morphism from a smooth
projective surface to a smooth non-rational curve. Assume that $-(dK_{S \vert C} + \Delta)$
is nef for some  rational number $d > 1$ and some effective reduced divisor $\Delta$ (possibly $0$).
Then $f$ is smooth and locally trivial, the general fiber $f$ has genus $g(F) \leq 1$ and 
one of the following cases occurs.
\sn (a) $g(F) = 1, \Delta = 0,$ and $-K_{S \vert C}$ is a torsion line bundle
\sn (b) $g(F) = 0,$ and every connected component $\Delta_i$ of $\Delta$ is a smooth curve, numerically
equivalent to $-rK_{S \vert C}$ for some positive rational number $r;$ moreover $\f : \Delta_i
\la C$ is \'etale.}

\bn Observe the following special case. If $f: S \la C$ is a $\bP_1-$bundle, then $-K_{S \vert C}$
is nef if and only if $S = \bP (E)$ with $E$ semi-stable or equivalently $E \otimes {{\rm det E^*}
\over {2}}$ is nef. We shall explain this in more detail and in any dimension in (4.6).   

\bn {\bf Proof.}  We shall proceed in several steps. From 
$$ 0\leq -\big( d \omega_{S/C} + \Delta\big) F \leq - d K_S F$$ we deduce that either 
$K_S F = 0$ and $F$ is elliptic, or $K_S F < 0$ and $F$ is rational. In the former case, we factor
 out $f$ as $\sigma \circ \pi$ where $\sigma : R \to C$ is a relatively minimal elliptic 
fibration and $\pi: S\to R$ is a birational morphism. Since $g(C) \ge 1$ we get $\kappa(R)\ge 0,$ 
and thus $\chi({\cal O}_R) \ge 0.$ The canonical bundle formula yields $$\omega_{R/C} = \sigma^*(D) 
+ \sum^t_{i=1}(m_i - 1) F_i$$ where $D$ is a divisor of degree equal to $\chi({\cal O}_R) \ge 0,$ 
and $m_1 F_1,\dots,m_t F_t$ stand for the multiple fibres of $\sigma$. Hence $\omega_{R/C}$ is nef. 
We also have $\omega_{S/C}= \pi^*\big(\omega_{R/C}\big) +E,$ for some effective divisor $E.$ If $L$ now 
stands for an ample divisor on $S$ we get $$0\le - \big(d \omega_{S/C} + \Delta\big) L = - d\, 
\pi^*\big(\omega_{R/C}\big) L - \big(d\, E + \Delta)\, L \le 0\,.$$
We conclude that $\Delta = 0,\ E = 0,\ f = \sigma$ and $\omega_{S/C}$ is numerically trivial. 
Proposition 4.2 applies here, and yields part~(i). 
\sn From now on we shall assume that $F$ is rational. \bigskip
\noindent {\bf Claim~1:} If $\beta : S \to R$ is the blow-down of a $(-1)$-curve $E,$ and 
$\Delta^\prime:= \beta (\Delta)$ is the reduced image of $\Delta,$ then $-(d\, \omega_{R/C} 
+ \Delta^\prime)$ is nef.

The proof is just a computation. We write 
$$\eqalignno{\Delta & = \beta^*(\Delta^\prime) - m\, E\,,\qquad m\ge 0\cr 
\omega_{S/C} & = \beta^*\big(\omega_{R/C}\big) + E\,.}$$ Then
$$-\big(d\, \omega_{S/C} + \Delta\big) = - \beta^* \big(d\, \omega_{R/C} + 
\Delta^\prime\big) + (m-d) E \leqno{(*)}$$ and $0\le - \big(d \omega_{S/C} + 
\Delta\big)\, E = d-m\,.$ 

Let $A^\prime$ be any irreducible curve in $R.$ Its strict transform is of the 
form $A = \beta^* (A^\prime) - r E\,,\ r\ge 0.$ Thus $$\eqalign{0\le - \big(d\, 
\omega_{S/C} + \Delta\big) A & = - \big(d\, \omega_{R/C} + \Delta^\prime\big) 
A^\prime + r (m-d)\cr &\le - \big(d\, \omega_{R/C} + \Delta^\prime\big) A^\prime\,.}$$ 
This finishes the proof of Claim (1).
\sn Now, let us assume for a moment that part~(b) of the Proposition holds true for smooth maps $f.$ 
Then, we are going to show that our $f$ is actually smooth. Otherwise, some $(-1)$-curve on 
$S$ could be blown down to a point $P\in R$ by a map $\beta: S\to R.$ From Claim~1 we know that 
$-\big(d\, \omega_{R/C} + \Delta^\prime\big)$ is nef. Since by a finite sequence of blow-downs 
we eventually reach a ${\bP}_1$-bundle, we may already assume that $R \to C$ is a smooth map. 
In this case, since we are assuming that the Proposition is true for $R,$ it follows that the 
multiplicity of $\Delta^\prime$ at $P$ is $m = 0$ or 1, and so $ (m-d)^2 > 0$. In view of $(*)$
 above we get
$$0\le \big(d\, \omega_{S/C} + \Delta\big)^2 = \big(d\, \omega_{R/C} + 
\Delta^\prime\big)^2 - (m - d)^2\,.\leqno(**)$$
Our assumptions again imply that $d\, \omega_{R/C} + \Delta^\prime$ is numerically equivalent to
$b\, \omega_{R/C},$ for some $b\in {\bQ}$. Then, $-\omega_{R/C}$ being nef combined with 
Propositon~4.1, yields $\omega^2_{R/C}=0,$ and so $$\big(d\, \omega_{R/C} + \Delta^\prime\big)^2 =0.$$ 
This is in contradiction to $(**).$ 

\sn It only remains to prove the Proposition in the particular case when $f$ is a ${\bP}_1$-bundle. 
Assume so in the sequel.
We shall freely use notation and results from~[Ha77], V.~2. Let $C_0$ stand for a section of $f: S\la C$ 
with minimal self-intersection $C^2_0 = -e.$ The decomposition of $\Delta$ into irreducible components
$$\Delta = C_1 + \dots+ C_r + F_1 + \dots + F_s + \sum_i \Delta_i$$ is written in such a way that 
the $C'_i s$ are exactly the components numerically equivalent to $C_0, F_1,\dots, F_s$ are the fibres 
in $\Delta$ and the remaining components are $\Delta_i \equiv a_i C_0 + b_i F\,,\ a_i, b_i$ integers. 
Since $\omega_{S/C}\equiv -2 C_0 - e F$ we get 
$$d\omega_{S/C} + \Delta \equiv (r-2d) C_0 + (s-de) F + 
\sum_i (a_i C_0 + b_i F)\,.$$ 
From $0\le \big(-d
\omega_{S/C} - \Delta\big) C_0$ it follows $$0\le (r-d) e-s + \sum_i (a_i e- b_i)\,.\leqno(***)$$ Let 
us consider the case $e > 0$ first. Since $C_0^2 = -e < 0$ we get $r = 0$ or 1, and so $(r-d) e < 0.$ 
Furthermore $a_i e-b_i \le 0$ for all $i$ ([Ha77,V.2]), which contradicts $(***)$. When $e=0$ we have $b_i \ge 0,$ 
and from $(***)$ it follows $$0\le -s- \sum_i\, b_i \le 0.$$ Hence $s$ and all $b'_i s$ are 0.
In particular, all components of $\Delta$ are numerically equivalent to a multiple of $C_0,$ and so is 
$\omega^{-1}_{S/C} \equiv 2 C_0,$ whence the claim.
\sn We shall finally deal with the case $e < 0.$ Now 
$b_i \ge a_i e/2$ for all $i$. Since $$- \omega_{S/C}\equiv 2 C_0 + e\, F, - \omega_{S/C}$$ is a limit of 
ample ${\bQ}$-divisors ([Ha77,V.2]), and thus it is nef. Note that $\omega^2_{S/C}=0$. Therefore
$$\eqalign{0&\le - \big(d \omega_{S/C} + \Delta\big) \big(-\omega_{S/C}\big) = \Delta \omega_{S/C} =\cr
&= r e - 2 s + \sum_i \big(a_i e - 2 b_i\big) \le 0\,.}$$ We conclude $r=s=0,\ a_i e-2 b_i = 0$ for all 
$i$, which yields the result. 

As for the fact that any component of $\Delta$ is mapping onto $C$ without ramification, it just follows 
from Hurwitz formula, namely $$2g(\Delta_i) -2 = \big(2g(C)-2\big) \Delta_i F = deg\, f^*(K_C)_{\vert\Delta_i} 
\cdot deg (\Delta_i \to C).$$

\bn {\bf Remark.} Let $X$ be a terminal variety of dimension $n,$ $Y$ a projective normal
${\bf Q}-$ Gorenstein variety
of dimension $n-1$ and $\f: X \la Y$ a surjective map such that $\XY$ is nef. Then 
$$ \dim \{y \in Y \vert X_y {\rm \ is \ singular}  \} \leq n-2.$$
In fact, this follows from (4.3) by taking $n-3$ general hyperplane sections.      

\bn We are now going to study threefolds $X$ admitting a map $\pi: X \la C$ to a curve
of genus at least 1 such that $\XC$ is nef. We start with a general statement, valid in
every dimension and generalising [Mi87,3.1].

\bn {\bf 4.5 Proposition} \ ÊÊ{\it Let $X$ be a n-dimensional projective manifold, $\pi: X \la C$
an extremal contraction to the smooth curve $C.$ Let $\lambda$ be the class of $\XC$ in
$N^1(X).$ Then the following statements are equivalent.
{\item {(1)} $\XC$ is nef
\item {(2)} the ample cone ${\overline {NA}}(X)$ is generated (as cone) by $\lambda $ and a
fiber $F$ of $\pi,$ i.e. $$\overline {NA}(X) = \bR_+ \lambda + \bR_+ F $$
\item {(3)} $\overline {NE}(X) = \bR_+ \lambda^{n-1} + \bR_+ \lambda^{n-2}ÊF $
\item {(4)} $\lambda^n \geq 0$ and every effective divisor in $X$ is nef.}}

\bn {\bf Proof.} First note that $\lambda$ and $F$ are clearly linearly independent in $N^1(X)$
and that moreover $\lambda^{n-1}$ and $\lambda^{n-2}F$ are linearly independent in $N_1(X).$ 
This statement holds because $\lambda^{n-1}F = (-K_F)^{n-1} > 0,$ whereas $\lambda^{n-2}F^2 = 0.$
\sn $(1) \Longrightarrow (2)$  This is clear since $\rho (X) = 2$ and $\lambda$ is nef
but not ample by (4.1).
\sn $(2) \Longrightarrow (3)$ One inclusion being obvious, we take an irreducible curve 
$C \subset X.$ Write in $N_1(X)$ :
$$ C = a \lambda^{n-1} + b \lambda^{n-2}F.$$
From $\lambda \cdot C \geq 0$ and $F \cdot C \geq 0$ we deduce via $\lambda^n =0$ (4.2) and
$\lambda^{n-1}F > 0$ that $a \geq 0, b \geq 0,$ from which our claim follows. 
\sn Since both (2) and (3) clearly imply (1), all three statements are equivalent.
\sn $(3) \Longrightarrow (4)$ By (1) and (4.2) we have $\lambda^n = 0.$  Let $D$ be an effective divisor. By (3) it is sufficient to
show 
{\item {(a)}Ê$D \cdot \lambda^{n-1}Ê\geq 0,$
\item {(b)}Ê$D \cdot \lambda^{n-2}F \geq 0.$}
(a) is clear, since $\lambda$ is nef by (1). (b) holds because $\lambda^{n-2}\cdot D \in
\overline {NE}(X),$ (again since $\lambda$ is nef) and since $F$ is nef.
\sn $(4) \Longrightarrow (1)$ We will show
$$ h^0(\Om) \geq 0$$
for large $m.$ By Riemann-Roch we have
$$ \chi(\Om) = {m^nÊ\over {n!}} \lambda^n + {m^{n-1} \over {(n-1)!}} \lambda^{n-1}ÊF + O(n-2).$$  
This can be reformulated as follows
$$
\chi \big(\omega_{X/C}^{-m}\big)
= {m^n\over n !} \lambda^n +
{m^{n-1}\over 2(n-1) !}
\lambda^{n-1} \big(- K_X\big) + O(n-2). $$
We also take into account
$$\lambda^{n-1} (-K_X) = \lambda^n + \big( 2g(C)-2 \big) \lambda^{n-1}F.$$
Since $\omega^{-1}_{X/C}$ is
$\pi$-ample, we have
$R^j\pi_* (\omega_{X/C}^{-m}) = 0$ for
$j\geq 1,\, m >> 0,$ and thus
$H^q\big(\omega^{-m}_{X/C}\big)
= H^q\big(\pi_*
\omega_{X/C}^{-m}\big) = 0$
for $q \ge 2,\ m >> 0.$

\noindent We therefore get:
\itemitem{(a)} If
$\lambda^n > 0,$
then $h^0\big(\omega^{-m}_{X/C}\big)
>0,\, m >> 0.$ Hence
$\lambda$ is nef, so that
$\lambda^n=0$ by (4.1), a
contradiction.
\itemitem{(b)} If $\lambda^n=0,$ we can conclude as before if $g(C) \geq 2.$ If $C$ is elliptic,
nothing can be concluded. We set $b=\inf  \big \{ \beta \in \bR \vert
\lambda + \beta F\ \hbox{is nef} \big \}$. If $b \leq 0,$ then $\lambda$ is nef, so assume $b > 0.$ 
Write $L = \lambda + b F$ then. If $L^n > 0$ then
$(L - \varepsilon F)^n > 0$
for $\varepsilon > 0$ small enough, so
that $r(L-\varepsilon F)$
is effective if $r >>0,$ hence nef against
the choice of~$L.$ Hence $L^n
=0.$ On the  other hand $L^n = \lambda^n + nb \lambda^{n-1}F,$ so that $b = 0,$ contradiction.

\bn {\bf (4.6) Remark} We expect that the condition $\lambda^n \geq 0$ in (4.5(4)) can be omitted.
In case $X$ is a $\bP_{n-1}$bundle over $C,$ this is easily verified as follows. Write
$X = \bP (E)$ with a rank $n-1-$bundle $E$ on $C.$ Then $\lambda$ is nef if and only if $E$
is semi-stable. Now we prove that in case $E$ is instable, then not every effective
divisor in $X$ is nef. Normalise $E$ such that $H^0(E) \ne 0$ but $H^0(E \otimes L) = $
for every line bundle $L$ on $C$ of negative degree. Since $E$ is instable, $E$ is not nef.
Let $s \in H^0(E), s \ne 0.$ Let $D \subset X$ be the associated divisor in $\cO_{\bP (E)}(1).$
Then $D$ is not nef, since $\cO_X(D) \simeq \cO_{\bP (E)}(1),$ which is not nef.

\bn {\bf 4.7 Proposition} {\it Let $X$ be a smooth projective threefold, $\pi: X \la C$
a surjective map to the curve $C$ of genus $g \geq 1.$ Assume $\omega^{-1}_{X \vert C}$ to be nef. Assume
furthermore that there exists a conic bundle $\f: X \la S$ and a map $f: S \la C$
such that $\pi = f \circ \f.$ Then $\pi $ is smooth. } 
 
\noindent {\bf Proof.} Let $\Delta \subset S$ denote the discriminant locus of $\varphi.$ 
For any curve $B\subset S$ we known 
that $$ \omega^{2}_{X/C}\, \cdot\, \varphi^{-1} (B) = - (4 \omega_{S/C} + \Delta)\  B$$ ([Mi83], p.~96), 
so that $- (4\omega_{S/C} + \Delta)$ is nef. 

If $\Delta = 0$ then $\varphi$ is a $\bP_1$-bundle and $\omega^{-1}_{S/C}$ is nef. 
Hence $f$ is smooth and $\pi$, being a composite of smooth maps, is also smooth. 
Suppose $\Delta\not= 0.$ By (4.4), $\Delta$ is a (possibly reducible) 
smooth curve, 
all of whose components map surjectively onto $C.$ The morphism $\pi$ can only fail to be smooth at points 
lying on $\varphi^{-1}(\Delta).$ In order to see that this will never happen, take any point 
$P\in \Delta.$ Since $\Delta$ is smooth at $P, \varphi^{-1}(P)$ is a pair of distinct 
lines meeting at $Q$ ([Be77], 1.2). $\pi$ is smooth at every point of $\varphi^{-1} (P)$ 
different from $Q$. Let us see that $\pi$ is smooth at $Q$ too. We take a small 
analytic neighbourhood $U \subset S$ of $P$ with local parameters $(s,t)$ such 
that $P=(0,0), \Delta$ is locally defined by $s=0$ and $f$ becomes the projection 
$(s,t) \la s.$ We can consider $\varphi^{-1}(U)$ as the hypersurface in 
$U \times{\bP}^2$ given by an equation
$$\sum_{0\leq i\leq j\leq 2} A_{ij} (s,t) X_i X_j = 0 \leqno(*)$$ where 
$(X_0: X_1: X_2)$ are the homogeneous coordinates of ${\bP}_2,$ and 
the $A_{ij}'s$ are analytic functions (see [Be77]). We can also arrange things such that $\varphi^{-1}(P)$ 
is given by the equation $X^2_1 + X^2_2=0,$ so that $Q$ is (1:0:0) in 
$\{P\} \times {\bP}^2.$ We introduce affine coordinates $x_1 = X_1/X_0,\ x_2 = X_2/X_0$ 
and transform $(*)$~into $$
A_{00} (s,t) + A_{01} (s,t) x_1 + A_{02}(s,t) x_2 + \sum_{1\leq i\leq j\leq 2} A_{ij} (s,t) x_i x_j = 0,\leqno(**) $$
Since $(**)$ becomes $x^2_1 + x^2_2 = 0$ for $s=t=0,$ we obtain 
$$A_{11}(0,0) = A_{22} (0,0) = 1\,,\quad \hbox{and}\quad A_{ij} (0,0)=0\quad \hbox{otherwise}\,.$$
The series expansion of $A_{ij}(s,t)$ around $(0,0)$ thus becomes for 
$ i= 1, 2: A_{ii} (s,t) = 1 + \big(a_{ii} s + b_{ii}t\big) + $ (terms of degree $\ge 2$ in $s,t$), 
otherwise: $A_{ij} (s,t) = a_{ij} s + b_{ij} t + $ (terms of degree $\ge 2$ in $s,t$). 
\sn Since $\Delta$ is the discriminant locus of $\varphi$ we get that $\det A_{ij} (s,t) = 0$ if 
and only if $s=0.$ Therefore, $\det A_{ij} (s,t) = s \cdot F(s,t)$ for some analytic function 
$F.$ Since the linear term of $\det A_{ij} (s,t)$ is $a_{00} s + b_{00}t$ we deduce $b_{00} = 0.$ 
On the other hand, the linear term of $(**)$ in all four variables $s, t, x_1, x_2$ is $a_{00}s + b_{00}t$. 
Hence, $X$ being smooth at $Q$ implies $a_{00}\not= 0.$ 
\sn Finally, the fibre of $\pi$ over $s=0$ is $$ G (s, x_1, x_2)
= A_{00}(0,t) + A_{01}(0,t) x_1 + A_{02}(0,t) x_2 + \sum_{1\le i\le j\le 2} A_{ij} (0,t) x_i x_j = 0$$
Since
${\partial G\over \partial s} (Q) = a_{00} \not=0,$ we finally conclude that $\pi^{-1}(0)$ is non-singular 
at $Q,$ as claimed. 

\bn In general we have the following conjecture for the relative situation, some special
cases of which we shall prove. 

\bn {\bf 4.8 Conjecture}  Let $\pi: X \la C$ be a surjective morphism from the 
smooth projective threefold $X$ to the smooth curve $C$ of genus $\geq 1.$ Assume
that $\XC$ is nef and that the general fiber of $\pi$ has Kodaira dimension $- \infty.$
Then $\pi$ is a submersion. More precisely, there exists a sequence
$$ X = X_0 \buildrel \f_1 \over \la X_1 \buildrel \f_2 \over \la X_2 \la ... \buildrel \f_r
\over \la X_r \leqno (4.8.1)$$
of birational morphisms over $C,$ each $\f_i$ being the blow-up of a smooth curve 
$C_i \subset X_i$ which map without ramification to $C,$ such that all $\omega^{-1}_{X_i
\vert C}$ are nef and the resulting map $f: X_r \la C$ is 
\sn (1) either a smooth Mori fibration, the fibers being del Pezzo surfaces (so that in particular
$\rho (X_r) = 2)$ or 
\sn (2) $f$ factors as $X_r \buildrel h \over \la S \buildrel g \over \la C,$
with $h$ a (Mori) conic bundle and $g$ a $\bP_1-$bundle (hence $\rho(X_r) = 3.$) 

\bn In case (2), $-(4\omega_{S \vert C} + \Delta)$ is nef by (4.7), and the ramification 
$\Delta $ of $h$ is described in (4.4). Note that in case $g(C) = 1$
the conjecture is an immediate consequence of our Main Theorem (and its corollaries).
In case $\pi$ is the Albanese map, we have proved (4.8) in (3.4). 
It turns out that, after suitable finite etale base change, the structure in the
above conjecture can be made quite simple (cp. (3.5)):

\bn {\bf 4.9 Proposition}  {\it Assume Conjecture (4.8) holds. Then after a suitable
\'etale base change $B \la C$ the induced submersion 
$$\sigma : X' = X \times_C B \la B$$ has the following structure.
\sn There exists a sequence
$$ X' = X'_0 \buildrel \f'_1 \over \la X'_1 \buildrel \f'_2 \over \la X'_2 \la Ê\ldots
\buildrel  \f'_r  \over \la X'_r $$
with the same properties as in (4.8) and $f' : X'_r \la B$ belongs to one of the 
following cases.
\sn (1) Either $\rho(X'_r) = 2$ and $f'$ is a $\bP_1 \times \bP_1-$bundle or a
$\bP_2-$bundle (in the latter case $X'_r = \bP(E)$ with a semistable rank 3-bundle
$E$ over $B$) or
\sn (2) $\rho(X'_r) = 3$ and $f'$ factors as $$X'_r \buildrel {h'} \over \la S' 
\buildrel {g'} \over \la B$$ where both $h$ and $g$ are $\bP_1-$bundles and moreover
$S' = \bP(F)$ with $F$ a semistable rank 2-bundle on $B.$}

\bn The proof of (4.9) is again an application of (0.4) and just the same of corollary (3.5)
which is contained in the proof of the Main Theorem. For the semi-stability of the
bundles in question apply [Mi87,3.1].

\bn {\bf (4.10)} Let $X$ be a smooth projective threefold and $\pi: X \la C$ a
surjective morphism to the smooth curve $C$ of positive genus. In order to prove
Conjecture 4.8 we need to investigate birational extremal contractions $\f: X \la W.$ As
in (3.2) above and [DPS93,p.234]  we see that $\f$ is the blow-up of a smooth
curve $C_0 \subset W.$ Since $g(C) > 0,$ we have a factorisation $\pi = \sigma \circ
\f$ with a map $\sigma : W \la C.$ In this situation we can state

\bn {\bf 4.11 Proposition} {\it Assume $\omega^{-1}_{W \vert C}$ nef. Let $S$ be the 
general fiber of $\pi$ and assume either
$K_S^2 > 0$ or $\dim \sigma (C_0) = 1.$ Then $\sigma \vert C_0 : C_0 \la C$ is \'etale.}

\bn {\bf Proof.} Let $E$ denote the exceptional divisor of $\f$ and let $F$ be a
(general) fiber of $E \vert C_0.$ Set 
$$g = g(C_0),Ê\gamma = g(C), d = {\rm deg}(C_0 \la C).$$
Here $d = 0$ iff $\dim \sigma (C_0) = 0.$ Following [DPS93,p.234-235] we write
for numerical equivalence
$$ -K_{X \vert E} \equiv C_1 + bF , \ N_{E \vert X}Ê\equiv -C_1 + \mu F.$$  
Then
$$ \omega_{X \vert C} \vert E = -C_1 -(b+d(2\gamma -2))F.$$ 
We know by (4.1) that $\omega^3_{X \vert C}Ê= \omega_{X \vert C}^3 = 0.$ Thus
$$ 0 = (\f^*(\omega_{W \vert C}))^3 = (- E)^3 = -3 (\omega_{X \vert C})^2 \vert E + 
3 (\omega_{X \vert C} \vert E \cdot (E \vert E) - (E \vert E)^2 $$
$$ = e - 3b - 6d(\gamma -1) - \mu, \leqno (1)$$
where $e = -C_1^2.$
Moreover (**) of [DPS93,p.235] gives
$$ e - b - 2(g-1) + \mu = 0. \leqno (2)$$ 
Since $\XC$ is nef, we have
$$ \XC \cdot C_1 \geq 0, \  \XC^2 \cdot E \geq 0,$$
which translate into
$$ e - b - 2d(\gamma-1) \leq 0 \leqno (3) $$
and
$$ e - 2b - 4d(\gamma - 1) \leq 0. \leqno (4)$$ 
The nefness of $\XC$ also yields
$$ 0 \geq \omega_{X \vert C} \cdot C_0 = \f^*(\omega_{W \vert C}) \cdot C_1 = (\omega_{X \vert C} - E) \cdot C_1 =
-b - 2d(\gamma - 1) - \mu$$
and therefore
$$ b + \mu + 2d(\gamma - 1) \geq 0. \leqno (5)$$
Note that by (1) 
$$ 0 = e - 2b - 4d(\gamma - 1) - (b + \mu + 2d(\gamma - 1)) \leq e - 2b - 4d(\gamma - 1)
\leq 0.$$
The first and second inequality are due to (5) and (4), respectively. Hence (4) and
(5) are just equalities:
$$ e - 2b - 4d(\gamma - 1)  = 0 $$
$$ b + \mu + 2d(\gamma - 1) = 0.$$ 
\sn We first deal with the case $d > 0.$ Note that $g-1 \geq d(\gamma - 1). $ Adding up (1)
and (3) we get
$$ 0 = 2e - 4b - 6d(\gamma -1) - 2(g-1) \leq 2e - 4b- 8d(\gamma - 1).$$
Then
$$ b + 2d(\gamma - 1) \leq {1 \over 2} e.$$
On the other hand we obtain from $\XC \cdot C_0$ that
$$ b + 2d(\gamma - 1) \geq {1 \over 2}e $$
resp.$$ b + 2d(\gamma - 1) \geq, $$ 
if $e > 0.$ We thus conclude 
$$ b + 2d(\gamma - 1) = {1 \over 2}e \leq  \ {\rm and} \ g-1 = d(\gamma - 1).$$
This implies that $C_0 \la C$ is \'etale of degree $d.$
\sn Now suppose $d = 0.$ Then we have $K_S^2 > 0$ by assumption. Combining (2),(6)
and (7) we deduce $g = 1$ and $b = {e \over 2}.$ Then (3) yields $e \leq 0,$ thus
$e  = 0$ or $-1.$ But $-1 = e = 2b $ is absurd. Hence
$$ e = b = \mu = 0, \ g = 1. \leqno (8).$$ 
On the other hand $E$ is contained in some fiber $S$ of $\pi,$ hence, taking into
account (6), we derive
$$ (\omega^{-1}_{X \vert C} + S)^2 E = (\omega^{-1}_{X \vert C})^2 E = 0.$$
Now the nef divisor $\XC + S$ is also big thanks to the assumption $K_S^2 > 0.$ 
Furthermore
$$ (\XC + S) E^2 = (-K_X) E^2 = C_1 \cdot (-C_1) = e = 0,$$
in view of (8). Then the following proposition gives 
$$ 0 \equiv (\XC + S) \vert E = C_1 $$
which is absurd. This concludes the proof.

\bn {\bf 4.12 Proposition} \ ÊÊ{\it Let $X$ be a projective manifold of dimension $n.$
Let $D$ be a nef and big divisor on $X$ and $E$ a divisor with $D^{n-1} \cdot E = 0.$
Then $D^{n-2}Ê\cdot E^2 \leq 0,$ with equality holding if and only if 
$D^{n-2}Ê\cdot E \equiv 0.$}

\bn {\bf Proof. } [Lu90].

\bn We now turn to the case of mappings to surfaces. 

\bn {\bf 4.13 Conjecture}  Let $X$ be a smooth threefold, $S$ a smooth surface and $\pi: X \la S$
a surjective map with connected fibers. If $\omega^{-1}_{X \vert S}$ is nef, then $\pi$ is smooth.

\bn Note that if $\omega^{-1}_X$ is nef, then $\pi$ may very well be non-smooth, e.g. there are
Fano threefolds which are conic bundles with non zero discriminant ove $\bP_2.$ But observe
that "$-K_{X\vert S}$ nef" is a somehow stronger condition than the nefness of $-K_X$ if
$\kappa (S) = - \infty.$ 

\bn The general fiber of $\pi$ is either elliptic or a rational. In the former case the conjecture
follows from 4.2. So we shall assume from now on that it is rational. In case $S$ is abelian,
4.13 is our Main Theorem. If $C$ is a general 
hyperplane section of $S$ and $X_C = \pi^{-1}(C),$ then $\omega^{-1}_{X_C \vert {C}} =
\XC \vert X_C$ is nef, and therefore $\pi$ is smooth over $C.$
Hence $\pi$ can fail to be smooth only over finitely many points of $Y.$
\sn The following is a straightforward 
consequence of (4.13).

\bn {\bf 4.14 Proposition}  {\it Let $\pi : X \la Y$ be a surjective morphism between 
projective manifolds with $\dim X = \dim Y + 1.$ Let $B \subset Y$ be the set points over which
$\pi$ fails to be smooth. Let $\omega_{X \vert Y}^{-1}$ is nef. If Conjecture (4.13) holds,
then ${\rm codim}_YB \geq 3.$}

\bn {\bf 4.15 Proposition} \ ÊÊ{\it In order to prove Conjecture 4.13, we may assume that $S$
contains no rational curve and that $H^1(S,\cO_S) \ne 0.$}

\bn {\bf Proof.} Let $B \subset Y$ be the set of point over which $\pi$ is not smooth. Take a
Lefschetz pencil $\Lambda$ of hyperplane sections on $S$ such $B$ is disjoint from the base locus.
Take a sequence of blow-ups, say $\beta_1 : R_1 \la S$ to make the map associated to $\Lambda$
base point free. We obtain a map $f : R_1 \la \bP_1$ with reduced fibers. Choose any smooth
hyperplane section $C$ of $R_1$ of positive genus, not passing through the singular fibers of
$f,$ nor through any point of $\beta_1^{-1}(B).$ We arrange things that $C \la \bP_1$ is
unramified where $R_1 \la \bP_1$ is not smooth. Then
$$ R_2 = C \times_{\bP_1}ÊR_1 $$ is a smooth surface which is mapped onto $C,$ so that $R_2$
contains at most a finite number of rational curves. The next step will be to choose a  
hyperplane section $D$ and a smooth curve $H \in \vert nD \vert,$ which skips the singular
points of all rational curves in $R_2$ and also avoids all points lying over $B.$ 
Let $R_3 \la R_2$ be the $n-$cyclic cover determined by $H.$ The rational curves in $R_2$
become irrational when lifted to $R_3,$ since $n \gg 0.$ Hence $R_3$ contains no rational
curves. 
\sn Let $\beta_i : R_i \la R_{i-1}$ the canonical map, $X_i \buildrel \pi_{i}Ê\over \la 
R_i$ the base change with associated maps $\alpha_i : X_i \la X_{i-1}$. Here we denote $X = X_0$  
and $\pi = \pi_1.$ If $\beta_1$ is the blow-up of $S$ at  $B = \{P_1, ..., P_r \},$ and if the
$E_i$ are the corresponding exceptional divisors in $R_1,$ then $\alpha_1$ is the blow-up
of $X$ at $\pi^*(E_1), \ldots , \pi^*(E_r).$ Since
$$ K_{R_1}Ê= \beta_1^*(K_S) + \sum E_i,$$ and
$$ K_{X_1} = \alpha_1^*(K_X) + \sum \pi^*(E_i),$$
we get $\omega^{-1}_{X_1 \vert R_1}Ê= \alpha^*(\omega^{-1}_{X \vert S}),$ which is nef.
From the fact that $C \la \bP_1$ is branched away from the singular points of $R_1 \la \bP_1,$
we obtain 
$$ \omega_{R_2 \vert C}Ê= \beta_2^*(\omega_{R_1 \vert \bP_1})$$
and $$ \omega_{X_2 \vert C}Ê= \alpha^*_2(\omega_{X_1 \vert \bP_1}).$$
We conclude that $\omega^{-1}_{X_2 \vert R_2})Ê= \alpha^*(\omega^{-1}_{X_1 \vert R_1}$,
hence $\omega^{-1}_{X_2 \vert R_2}$ is nef. Now $\alpha_3$ is a $n-$cyclic cover totally
ramified at $\pi^*(H)$ and determined by $\pi^*(H) \sim n\pi^*(D).$ From e.g. [BPV84,p.42]
we derive
$$ K_{R_3}Ê= \beta_3^*(K_{R_2} + (n-1)D) $$
and
$$ K_{X_3}Ê= \alpha^*_3(K_{X_2}Ê+ (n-1)\pi^*(D)).$$
Hence $\omega^{-1}_{X_3 \vert R_3}Ê= \alpha^*_3(\omega^{-1}_{X_2 \vert R_2})$ is nef. 
By construction the map $X_3 \la X$ is \'etale over the singular fibers of $\pi.$ 
If therefore we can show that $\pi_3$ is smooth, then $\pi$ is smooth, too.

\bn {\bf (4.16)} In view of the preceding result, we can restrict ourselves to the 
following situation. $X$ is a smooth projective threefold, $S$ a smooth surface without
rational curves and such that $q(S) > 0.$ Let $\pi: X \la S$ be surjective with connected
fibers. Assume that the general fiber is rational.
\sn Since $K_X$ is not nef, there exists an extremal contraction $\f : X \la W.$ We are
going to investigate the structure of $\pi.$

\bn {\bf 4.17 Proposition}  {\it In the situation of (4.16) assume $\dim W \leq 2.$ Then 
$W = S, \f = \pi $ and $\pi$
is a $\bP_1-$bundle. }

\bn {\bf Proof.} Since $S$ does not contain rational curves, it is clear that $\dim W = 2$ and that
there is a map $\sigma : W \la S$ such that $\pi = \sigma \circ \f.$ Since the fibers of
$\pi$ are connected, $\sigma$ must be birational, i.e. a sequence of blow-ups. Let $E$ be
the exceptional divisor of $\sigma$ and $\Delta$ the discriminant locus of the conic bundle
$\f.$ Then an easy calculation shows (cp. 1.6, 1.7)
$$ 0 \leq \omega^{-1}_{X \vert S} \cdotÊ\f^*(C) = - (\Delta + 4E) C.$$
Hence $\Delta = E = 0$ and the claim follows.

\bn It remains to treat the case that $\f$ is birational. Let $E$ be the exceptional divisor.
\bn {\bf 4.18 Proposition}  {\it $\dim \f(E) = 1.$}

\bn {\bf Proof.} Assume $\dim \f(E) = 0.$ Similar as in [DPS93,3.3] we see 
that $\omega_{W \vert S}^{-1}$ is big and nef. This contradicts (4.1). 

\bn So $\f$ is the blow-up of a smooth curve $C_0.$ Since $S$ does not contain rational
curves, we obtain again a factorisation $ \pi = \sigma \circ \f,$ with $\sigma : W \la S.$
 
\bn {\bf 4.19 Proposition} \  {\it
{\item {(1)} $\dim \sigma (C_0) = 0.$
\item {(2)} $C_0 \simeq \bP_1.$}}

\bn {\bf Proof.} (1) follows easily from the
remarks after (4.13). 
\sn (2) Choose $H$ ample on $S$ and set $L = \omega^{-1}_{X \vert S} + \pi^*(H).$
Then $L$ is nef and big and
$$ aL - K_X = (a+1) \omega^{-1}_{X \vert S} + \pi^*(aH - K_S)$$
is also nef and big for $a \gg 0.$ Therefore $mL$ is generated by global sections for
large$m$ by the base point free theorem. Let $D \in \vert mL \vert$ be a general smooth and
irreducible element. Let $A = D \cap E.$ We may assume $A$ smooth and irreducible. Let $f = \pi
\vert D : D \la S.$ Then $f$ is generically finite and by (1) $A$ is contracted by $f.$
Therefore $(A^2)_D < 0.$ On the other hand,
$$ 0 > (A^2)_R = E^2 \cdot mL = m E \vert E \cdot L \vert E = m E \vert E \cdot (-K_X) \vert E.$$
In the notations of the proof of (4.11) we obtain the following inequality
$$ 0 > m(-C_1 + \mu F)(C_1 + bF) = m(e + \mu - b) = 2m(g - 1),$$
where $g$ is the genus of $C_0.$ Consequently $g = 0.$ 

\bn {\bf 4.20 Proposition} \ ÊÊ{\it Suppose we know the following 
\sn (*) Let $Z$ be a smooth projective threefold  having a surjective morphism $f : Z \la Y$
to a smooth surface $Y$ having no rational curve. Assume $q(S) > 0.$ and  that $-K_{Z \vert Y}$ nef.
If $g : Z \la Z' $ is a birational extremal contraction, then $-K_{Z' \vert Y}$ is again
nef. 
\sn Then in our situation $\pi : X \la S$ is a submersion and in particular $\f$ cannot exist. }

\bn {\bf Proof.} In view of (4.16) we have a birational contraction $\f : X \la W$  contracting
the divisor $E$ to the curve $C_0 \subset W.$ Moreover there $\omega^{-1}_{W \vert S}$ is nef
via the induced map $\sigma : W \la S.$ Again we shall use the notations of the proof of
(4.11). In the same way as (4.11(1)) we get
$$ e - 3b - \mu = 0. \leqno (1)$$ 
Since $C_0$ is rational, (**) of [DPS93,p.235] gives
$$ e - b + \mu = -2. \leqno (2)$$ 
Adding up (1) and (2) gives
$$ e - 2b = -1. \leqno (3)$$ 
Since $\omega_{X \vert C} \vert E = (-K_X) \vert E = C_1  + bF$ is nef, we obtain $b \geq  e.$
Combining with (3) yields $e = 1,$ hence $b= 1, \mu = -2.$ 
In view of our hypothesis we can apply this procedure inductively finitely many times
until we reach the situation where no birational contraction is possible on $W.$ 
From (4.17) it follows that $\sigma : W \la S$ is a $\bP_1-$bundle. Since then $C_0$ is
contracted to a point by $\sigma,$ it is a fiber of $\sigma$ and thus $N_{C_0 \vert W} = 
\cO \oplus \cO.$ This contradicts $ e = 1.$

\bn The condition (*) is "mostly" satisfied as we explain in the next two propositions which
are proved with the same type of arguments
as Propositions~(3.3)
and~(3.5) in~[DPS93], respectively:
\bigskip
\noindent {\bf 4.21 Proposition} {\it Let $X$
be a smooth projective threefold and let
$\pi: X \rightarrow Y$ be
a surjective morphism with connected
fibres, where $Y$ is either a smooth
curve of genus
$\ge 1$ or a smooth irregular surface
containing no rational curve. Suppose
the general
fibre of $\pi$ has Kodaira
dimension $-\infty$. Let
$\varphi: X \rightarrow W$ be the
blow-up of a smooth curve $C_0 \subseteq W.$ 
We always have a factorization
$\pi= \sigma \circ
\varphi,$ where
$\sigma : W \rightarrow C.$ Now assume that $\omega^{-1}_{W/Y}$ is not nef.
Then $C_0 \simeq \bP_1, \sigma(C_0)$
is a point and one of the
following cases occurs:
\smallskip
\item{(A)} $N_{C_0/W} = \cO(-2)
\oplus \cO(-2),$ and $K_W \cdot C_0 = 2$
\item{(B)} $N_{C_0/W} = \cO(-1)
\oplus \cO(-2),$
and $K_W \cdot C_0 =1.$}
\bigskip
\noindent {\bf 4.22 Proposition}  {\it Case $B$ above is 
impossible.}
\medskip
\noindent
\noindent{\bf Proof} {\it of (4.22).}
We proceed exactly as in Proposition~(3.5)
in~[DPS93], replacing everywhere
$K_X, K_W$ by $\omega_{X/Y}, \omega_{W/Y},$ etc. At the end we obtain a threefold $Z$
with one terminal singularity such that the 
$\bQ$-divisor $\omega^{-1}_{Z/Y}$
is big and nef. Now we apply
[KMM87, 1.2.5, 1.2.6] 
with $\Delta = 0, D = f^*K_Y,$ where
$f: Z\rightarrow Y.$ It follows that
$H^1(Z, f^* K_Y) =
0$. The Leray spectral sequence yields
$H^1(Y, K_Y) = 0$,
which contradicts our hypothesis.

\vfill \eject 
{\medium References}

\bn 
{\item {[BPV84]}ÊBarth,W.;Peters,C.;Van de Ven,A.: Compact complex surfaces. Springer 1984
 
\item {[Be77]} Beauville,A.: Vari\'et\'es de Prym et jacobiennes interm\'ediaires. 
Ann. Sci. Norm. Sup. 10, 309-391 (1977)

\item {[Be83]}ÊBeauville,A.: Vari\'et\'es k\"ahl\'eriennes dont la premi\`ere classe de Chern
ets nulle. J. Diff. Geom. 18, 755-782 (1983)

\item {[CP91]} Campana,F.;Peternell,Th.: Projective manifolds whose tangent bundles are 
numerically effective. Math. Ann. 289, 169-187 (1991)

\item {[Cu88]} Cutkosky,S.: Elementary contractions of Gorenstein threefolds. Math. Ann. 280,
521-525 (1988)

\item {[DPS93]} Demailly,J.P.;Peternell,Th.;Schneider,M.: K\"ahler manifolds with numerically
effective Ricci class. Comp. math. 89, 217-240 (1993) 

\item {[DPS94]} Demailly,J.P.;Peternell,Th.;Schneider,M.: Compact complex manifolds with
numerically effective tangent bundles. J. Alg. Geom. 3, 295-345 (1994)

\item {[DPS96]} Demailly,J.P.;Peternell,Th.;Schneider,M.: Compact K\"ahler manifolds with
hermitian semipositive anticanonical bundle. Comp.math. 101, 217-224 (1996)

\item {[DPS97]} Demailly,J.P.;Peternell,Th.;Schneider,M.: In preparation

\item {[El82]} Elencwajg,G.: The Brauer group in complex geometry. Lecture Notes in
Math. 917, 222-230 (1982)

\item {[Fl87]} Fletcher,A.R.: Contributions to Riemann-Roch on projective 3-folds with
only canonical singularities and applications. Proc. Symp. Pure Math. 46, 221-231 (1987)

\item {[Fu78]} Fujita,T.: On K\"ahler fiber spaces over curves. J.Math.Soc. Japan 30, 779-794
(1978)

\item {[GR70]} Grauert,H.;Riemenschneider,O.: {\obeylines Verschwindungss\"atze f\"ur analytische Kohomologie- 
gruppen auf komplexen R\"aumen. Inv. math. 11, 263-292 (1970) 

\item {[Ha77]} Hartshorne,R.: Algebraic geometry. Graduate texts in mathematics. Springer 1977

\item {[Ka81]} Kawamata,Y.: Characterisation of abelian varieties. Comp. math. 43, 253-276
(1981)

\item {[Ka85]} Kawamata,Y.: Minimal models and the Kodaira dimension of algebraic fiber spaces.
J. reine u. angew. Math. 363, 1-46 (1985)

\item {[KMM87]} Kawamata,Y.;Matsuda,K.;Matsuki,K.: Introduction to the minimal model problem.
Adv. Stud. Pure Math. 10, 283-360 (1997)

\item {[Ka96]} Kawamata,Y.: {\obeylines Divisorial contractions to 3-dimensional terminal quotient singularities.
Proc. Higher Dimensional Complex Geometry Trento 1994. Ed. M.Andreatta,
Th.Peternell. de Gruyter 1996}

\item {[KoMiMo92]} Koll\'ar,J.;Miyaoka,Y.;Mori,S.: Rationally connected varieties. J. Alg. Geom. 1,
429-448 (1992)

\item {[Ko89]} Koll\'ar,J.: Flips, flops and minimal models. Surv. in Diff. Geom. 1, 113-199
(1991)

\item {[Mi83]} Miyanishi,M.: Algebraic methods in the theory of algebraic threefolds. Adv. Stud. Pure
Math. 1, 69-99 (1983)

\item {[Mi87]} Miyaoka,Y.: The Chern classes and Kodaira dimension of a minimal variety. Adv. Stud. Pure 
Math. 10, 449-476 (1987) 

\item {[Mi93]} Miyaoka,Y.: Relative deformations of morphisms and application to fibre spaces.
Comm. Math. Univ. Sancti Pauli 42, 1-7 (1993)

\item {[Mo82]}ÊMori,S. : Threefolds whose canonical bundles are not numerically effective. Ann. Math
116, 133-176 (1982)

\item {[Mo87]} Mori,S.: Classification of higher-dimensional varieties. Proc. Symp. 
Pure Math. 46, 269-331 (1987)

\item {[MP97]}ÊMiyaoka,Y.;Peternell,Th.: Geometry of higher dimensional varieties. DMV Seminar vol. 26.
Birkh\"auser 1997 

\item {[Lu90]} Luo,T.: A note on the Hodge index theorem. manusr. math. 67, 17-20 (1990)

\item {[Pe93]} Peternell,Th.: Minimal varieties with trivial canonical classes, I. Math. Z. 217, 
377-407 (1994)

\item  {[Re87]} Reid,M.: Young person's guide to canonical singularities. Proc. Symp.
Pure Math. 46, 345-414 (1987)

\bn \bn \bn 
{\obeylines
Thomas Peternell
Mathematisches Institut 
Universit\"at Bayreuth
D-95440 Bayreuth 
Germany
email: thomas.peternell@uni-bayreuth.de

\end

\end